\documentclass[11pt,oneside,onecolumn,reqno]{amsart}

\usepackage{amssymb,latexsym,mathrsfs,
amsxtra,amscd,amsfonts,amsthm,amsmath,verbatim,epsfig}

\usepackage{color}

\usepackage[colorlinks,pagebackref=true]{hyperref}

\usepackage[english]{babel}

\makeatletter

\makeatletter
\def\widebreve#1{\mathop{\vbox{\m@th\ialign{##\crcr\noalign{\kern3\p@}%
      \brevefill\crcr\noalign{\kern3\p@\nointerlineskip}%
      $\hfil\displaystyle{#1}\hfil$\crcr}}}\limits}

\def\brevefill{$\m@th \setbox\z@\hbox{$\braceld$}%
  \bracelu\leaders\vrule \@height\ht\z@ \@depth\z@\hfill\braceru$}
\makeatletter

\def\@citecolor{blue}
\def\@linkcolor{blue}
\def\@urlcolor{blue}
\def\@urlcolor{blue}

\numberwithin{equation}{section}
\def\im{\operatorname{im}}

\def\ker{\operatorname{ker}}
\def\dim{\operatorname{dim}}
\def\depth{\operatorname{depth}}

\def\ZZ{\mathbb Z}

\newcommand{\NN}{\mathbb N}
\newcommand{\mm}{\mathfrak m}

\newcommand{\ov}{\overline}
\newcommand{\lm}{{\lambda}}

\newcommand{\R}{\mathcal R}

\newcommand{\I}{\underline I}
\newcommand{\n}{\underline{n}}

\newcommand{\PP}{\mathcal P}

\newcommand{\kr}{\underline {k}}
\newcommand{\m}{\underline {m}}
\newcommand{\rrr}{\underline{r}}
\newcommand{\crr}{complete reduction }
\newcommand{\gcrr}{good complete reduction }
\newcommand{\ii}{i=1,\ldots,s}
\newcommand{\idl}{{I_1},\ldots,{I_s}}
\newcommand{\id}{{I_1}\cdots{I_s}}

\newcommand{\lmc}{C.(\underline{y}^{[l]},\fil(\n))}

\newcommand{\lf}{\left(}
\newcommand{\rg}{\right)}
\newcommand{\fil}{\mathcal F}
\newcommand{\gi}{{G_i(\fil)}_{++}}
\newcommand{\ggi}{{G_i(\fil)}}

\newcommand{\yy}{{y_1},\ldots,{y_d}}
\newcommand{\yk}{{y_1},\ldots,{y_k}}

\newcommand{\po}{P_{\fil}}
\newcommand{\ho}{H_{\fil}}
\newcommand{\hf}{H_{\I}}

\newcommand{\pf}{P_{\I}}
\newcommand{\bl}{\begin{lemma}}
\newcommand{\el}{\end{lemma}}
\newcommand{\bt}{\begin{theorem}}
\newcommand{\et}{\end{theorem}}
\newcommand{\ben}{\begin{enumerate}}
\newcommand{\een}{\end{enumerate}}
\newcommand{\bpf}{\begin{proof}}
\newcommand{\eepf}{\end{proof}}
\newcommand{\beqn}{\begin{eqnarray*}}
\newcommand{\eeqn}{\end{eqnarray*}}
\newcommand{\beqnn}{\begin{eqnarray}}
\newcommand{\eeqnn}{\end{eqnarray}}
\newcommand{\bd}{\begin{definition}}
\newcommand{\ed}{\end{definition}}
\newcommand{\bp}{\begin{proposition}}
\newcommand{\ep}{\end{proposition}}
\newcommand{\bc}{\begin{corollary}}
\newcommand{\ec}{\end{corollary}}
\newcommand{\bex}{\begin{example}}
\newcommand{\eex}{\end{example}}

\newcommand{\wrt}{with respect to }
\newcommand{\CM}{Cohen-Macaulay }
\newcommand{\wlg}{without loss of generality }

\theoremstyle{plain}
\newtheorem{theorem}{Theorem}[section]
\newtheorem{corollary}[theorem]{Corollary}
\newtheorem{proposition}[theorem]{Proposition}
\newtheorem{lemma}[theorem]{Lemma}

\newtheorem{example}[theorem]{Example}
\newtheorem{definition}[theorem]{Definition}

\theoremstyle{remark}
\newtheorem{remark}[theorem]{Remark}
\numberwithin{equation}{theorem}

\begin{document}

\title[]{Postulation and Reduction Vectors of Multigraded Filtrations of ideals}

\author{Parangama Sarkar}
\address{Department of Mathematics, Indian Institute of Technology Bombay, Mumbai, 400076, India}
\email{parangama@math.iitb.ac.in}
\author{J. K. Verma}
\address{J. K. Verma, Department of Mathematics, Indian Institute
of Technology Bombay, Mumbai, India}
\email{jkv@math.iitb.ac.in}

\begin{abstract} 
We study relationship between postulation and reduction vectors of admissible multigraded filtrations $\mathcal F= \{\fil (\n)\}_{\n\in\ZZ^s}$ of ideals in Cohen-Macaulay local rings of dimension at most two. This is enabled by a suitable generalisation of the Kirby-Mehran Complex. An analysis of its homology leads to an analogue of Huneke's Fundamental Lemma which plays a crucial role in our investigations.
We also clarify the relationship between the Cohen-Macaulay property of the multigraded Rees algebra of $\fil$ and reduction vectors with respect to complete reductions of $\fil.$ 
\end{abstract}

\thanks{
The first author is supported by CSIR Fellowship of Government of India.}
\maketitle
\thispagestyle{empty}

\section{Introduction}
The objective of this paper is to study properties  of Hilbert functions and Hilbert polynomials of multigraded filtrations of ideals under certain cohomological conditions. Among the themes presented
are: (1) an analogue of Huneke's Fundamental Lemma in terms of the homology of generalised Kirby-Mehran complex for multigraded filtrations of ideals using complete reductions (2) relationship between postulation vectors and reduction vectors for multigraded filtrations of ideals in Cohen-Macaulay local rings of dimension at most two (3) providing necessary and sufficient conditions for the equality of multigraded Hilbert functions and polynomials in terms of reduction numbers with respect to complete reductions and finally (4) relationship between the Cohen-Macaulay property of the Rees algebra of multigraded filtrations of ideals and reduction numbers in two-dimensional Cohen-Macaulay local rings.
\\Hilbert functions of mutigraded filtrations of ideals have been found useful in works of
B. Teissier \cite{T} who used them  in his investigations of Milnor numbers of singularities of complex analytic hypersurfaces. To wit, let $(R,\mm)$ be a $d$-dimensional local ring of an isolated singularity  of a complex analytic  hypersurface and let $f$ be the defining polynomial. Then the Jacobian ideal $J:=J(f)$ is an $\mm$-primary ideal and the function $B(r,s)=\lm(R/\mm^rJ^s)$ is a polynomial $P(r,s)$ of degree $d$ in  $r$ and $s.$ Here $\lm$ denotes length. Teissier proved that the normalised coefficients of monomials of degree $d$ in $P(r,s)$ are the Milnor numbers of linear sections of the isolated singularity. Joint reductions and the Bhattacharya function $B(r,s)=\lm(R/I^rJ^s)$ for $\mm$-primary ideals $I$ and $J$ have been used by D. Rees in several contexts. For example, he characterised pseudo-rational local rings of dimension two in terms of the constant term of the normal Hilbert polynomial for the normal Hilbert function $\lm(R/\ov{I^rJ^s}).$
\\We now describe the contents of the paper. We recall a few definitions and set up notation to explain the results of this paper.
\\Throughout this paper let $(R,\mm)$ be a Noetherian local ring of dimension $d$ and $\idl$ be 
 $\mm$-primary ideals of $R.$ For $s\geq 1,$ we put $e=(1,\ldots,1),\; \underline{0}=(0,\ldots,0)\in{\ZZ}^s$ and $e_i=(0,\ldots,1,\ldots,0)\in{\ZZ}^s$ where $1$ occurs at $i$th position. Let $\n=(n_1,\ldots,n_s)\in{\ZZ}^s,$ then we write ${\I}^{\n}=I_{1}^{n_1}\cdots I_{s}^{n_s}$ and $\n^+=(n_1^+,\ldots,n_s^+)$ where 
 \[ n_i^+ = \left\{
  \begin{array}{l l}
   {n_i}  & \quad \text{if ${n_i}>0$ }\\
    0 & \quad \text{if ${n_i}\leq 0.$ }
  \end{array} \right. \] 
  For $\alpha=({\alpha}_1,\ldots,{\alpha}_s)\in{\NN}^s,$ we put $|\alpha|={\alpha}_1+\cdots+{\alpha}_s.$ Define $\m=(m_1,\ldots,m_s)\geq\n=(n_1,\ldots,n_s)$ if $m_i\geq n_i$ for all $\ii.$ By the phrase ``for all large $\n$" we mean $\n\in\NN^s$ and $n_i\gg 0$ for all $\ii.$
\bd A set of ideals $\fil=\lbrace\fil(\n)\rbrace_{\n\in \ZZ^s}$ is called a $\ZZ^s$-graded {\bf{$\I=(\idl)$-filtration}} if for all $\m,\n\in\ZZ^s,$
{\rm (i)} ${\I}^{\n}\subseteq\fil(\n),$
 {\rm (ii)} $\fil(\n)\fil(\m)\subseteq\fil(\n+\m)$  and {\rm (iii)} if $\m\geq\n,$ $\fil(\m)\subseteq\fil(\n)$.
\ed
Let $t_1,\ldots,t_s$ be indeterminates. For $\n\in\ZZ^s,$ put $\underline t^{\n}=t_{1}^{n_1}\cdots t_{s}^{n_s}$ and denote the $\NN^s$-graded {\bf{Rees ring of $\fil$}} by $\mathcal{R}(\fil)=\bigoplus\limits_{\n\in \NN^s}
{\fil}(\n){\underline{t}}^{\n}$ and the $\ZZ^s$-graded {\bf{extended Rees ring of $\fil$}} by 
$\mathcal{R}'(\fil)=\bigoplus\limits_{\n\in{\ZZ}^s}{\fil}(\n){\underline{t}}^{\n}.$ For an $\NN^s$-graded ring $S=\bigoplus\limits_{\n\geq\underline 0}S_{\n},$ we denote the ideal $\bigoplus\limits_{\n\geq e}S_{\n}$ by $S_{++}.$ For $\fil=\{\I^{\n}\}_{\n \in \ZZ^s}$, we set $\mathcal R(\fil)=\mathcal R(\I),$ $\mathcal R^\prime(\fil)=
\mathcal R^\prime(\I)$ and $\mathcal R(\I)_{++}=\R_{++}.$
\bd A $\ZZ^s$-graded $\I=(\idl)$-filtration $\fil=\lbrace\fil(\n)\rbrace_{\n\in \ZZ^s}$ of ideals in $R$ is 
called an $\I=(\idl)$-{\bf{admissible filtration}} if ${{\fil}(\n)}={\fil}(\n^+)$ for all $\n\in \ZZ^s$ and $\mathcal{R}'(\fil)$ is a finite $\mathcal{R}'(\I)$-module. \ed
For an $\mm$-primary ideal $I$, the {\bf{Hilbert function}} $H_I(n)$ is defined as $H_I(n)=\lm\lf{R}/{I^n}\rg$ for all $n\in\ZZ.$ Here we adopt the convention that $I^n=R$ if $n\leq 0.$ P. Samuel \cite{samuel} showed that for sufficiently large $n,$ $H_I(n)$ coincides with a polynomial $P_I(n)$ of degree $d,$ called the {\bf{Hilbert polynomial}} of $I.$ For all $n\in{\ZZ},$ $P_I(n)$ is often written in the form $$P_I(n)=\displaystyle\sum_{i=0}^{d}(-1)^{i}e_i(I)\binom{n+d-1-i}{d-i}.$$ The coefficients $e_i(I)$ are integers for all $i=0,1,\ldots,d,$ called the {\bf Hilbert coefficients} of $I.$ The leading coefficient $e_0(I)$ is sometimes denoted by $e(I)$ and called the multiplicity of $I.$
\\ Let $s\geq 2$ and $\fil=\lbrace\fil(\n)\rbrace_{\n\in \ZZ^s}$ be an $\I$-admissible filtration of ideals in a Noetherian local ring $(R,\mm)$ of dimension $d.$ For the {\bf{Hilbert function}} $\ho(\n)=\lm\lf{R}/{{\fil}(\n)}\rg$ of $\fil,$ Rees \cite{rees3} proved that there exists a polynomial of degree $d$, called the {\bf{Hilbert polynomial}} of $\fil,$
$$\po(\n)=\displaystyle\sum\limits_{\substack{\alpha=({\alpha}_1,\ldots,{\alpha}_s)\in{\NN}^s \\ 
|\alpha|\leq d}}(-1)^{d-|{\alpha}|}{e_{\alpha}}(\fil)\binom{{n_1}+{{\alpha}_1}-1}{{\alpha}_1}
\cdots\binom{{n_s}+{{\alpha}_s}-1}{{\alpha}_s}$$ such that 
$\po(\n)=\ho(\n)$ for all large $\n.$ Here $e_{\alpha}(\fil)$ are integers called the {\bf{Hilbert coefficients}} of $\fil.$ This was proved by P. B. Bhattacharya for the filtration $\fil=\{I^rJ^s\}_{r,s\in\ZZ}$ in \cite{B} where $I$ and $J$ are $\mm$-primary ideals. Teissier \cite{T} showed the existence of $\po(\n)$ for the filtration $\fil=\{\I^{\n}\}_{\n\in\ZZ^s}.$
\\Let $I$ be an $\mm$-primary ideal in a local ring $(R,\mm)$ of dimension $d\geq 1.$ An integer $n(I)$ is called the {\bf{postulation number of I}} if $P_I(n)=H_I(n)$ for all $n> n(I)$ and $P_I(n(I))\neq H_I(n(I)).$ An ideal $J \subseteq I$ is
called a {\bf reduction} of $I$ if $JI^{n} = I^{n+1}$ for some $n$. We
say $J$ is a {\bf minimal reduction} of $I$ if whenever $K \subseteq J$ and
$K$ is a reduction of $I$, then $K=J$. 
Let $${r_J}(I)=\min\lbrace m:JI^n=I^{n+1}\hspace{0.1cm}
\mbox{ for }\hspace{0.1cm} n\geq m\rbrace\mbox{ and }
{r}(I)=\min\lbrace {r_J}(I):J\mbox{ is a minimal reduction of } I\rbrace.$$
The analogues of reduction and postulation number for $\I$-admissible multigraded filtration $\fil=\lbrace\fil(\n)\rbrace_{\n\in \ZZ^s}$ are described below. 
\bd
A vector $\n\in{\ZZ}^s$ is called a {\bf{postulation vector of $\fil$}} if $\ho(\m)=\po(\m)$ for all $\m\geq\n.$
\ed
Rees \cite{rees3} introduced the concept of \crr for $\lbrace{{\I}^{\n}}\rbrace_{\n\in{\ZZ}^s}.$ In a similar way we define \crr of an $\I$-admissible filtration $\fil=\lbrace\fil(\n)\rbrace_{\n\in \ZZ^s}.$ 
 \bd
A set of elements 
$\mathcal{A}={\lbrace}{x_{ij}\in{I_i}: j=1,\ldots,d;\ii}\rbrace$ is called a {\bf{complete reduction}} of $\fil$ if $J=(\yy)$ satisfies $J{\fil(\n)}=\fil(\n+e)$ for all large $\n\in\NN^s$ where $y_j=x_{1j}\cdots x_{sj}$ for all $j=1,\ldots,d.$
\ed
\bd
 A complete reduction $\mathcal{A}={\lbrace}{x_{ij}\in{I_i}: j=1,\ldots,d;\ii}\rbrace$ of $\fil$ is called a {\bf{good complete reduction}} if for all large $\m\in\NN^s,$ ${{{\fil}}(\m)}\cap(y_{1})=y_{1}{{{\fil}}(\m-{e})}$ where $y_1=x_{11}\cdots x_{s1}.$\ed
Let $\mathcal {A}={\lbrace}{x_{ij}\in{I_i}: j=1,\ldots,d;\ii}\rbrace$ be a complete reduction $\fil,$ $y_j=x_{1j}\cdots x_{sj}$ for all $j=1,\ldots,d$ and $J=(y_1,\ldots,y_d).$ 
\bd
A vector $\rrr\in\NN^s$ is called a {\bf{reduction vector of $\fil$ with respect to $\mathcal{A}$}} if for all $\n\geq\rrr,$ $J{\fil(\n)}=\fil(\n+e).$\ed
\bd An integer $k\in\NN$ is called the {\bf{complete reduction number of $\fil$ \wrt $\mathcal {A}$}} if 
$J{\fil(\n)}=\fil(\n+e)\mbox{ for all } \n\geq ke$ and whenever $k\neq 0,$ there does not exist any $0\leq t<k$ such that $\mbox{ }J{\fil(\n)}=\fil(\n+e)\mbox{ for all }\n\geq te.$
\ed
We use the following notation
\ben
{
\item[(1)]$\mathcal{P}(\fil)=\lbrace\n\in{\ZZ}^s|\mbox{ }\n\mbox{ is a postulation vector of }\fil\rbrace.$ 
\item[(2)]${\mathcal {R}}_{\mathcal A}(\fil)=\lbrace\n\in{\ZZ}^s|\mbox{ }\n\mbox{ is a reduction vector of }\fil\mbox{ with respect to  }\mathcal {A}\rbrace.$
\item[(3)] ${r}_{\mathcal A}(\fil)=$ The complete reduction number of $\fil$ \wrt $\mathcal{A}.$ }
\een
We now describe the main results proved in this paper. In section two, we prove some preliminary results about the coefficients of Hilbert polynomial of a multigraded filtration of ideals which we use to prove our main results.
Let $f(\n):{\ZZ}^s\rightarrow {\ZZ}$ be an integer valued function. Define the {\bf{first difference function}} of $f(\n)$ by ${\Delta}^1(f(\n))=f(\n+e)-f(\n).$ For all $k\geq 2,$ we define ${\Delta}^k(f(\n))={\Delta}^{k-1}({\Delta}^{1}(f(\n))).$ In \cite{huneke}, C. Huneke proved the following fundamental lemma:
\bl\cite[Lemma 2.4]{huneke}
Let $(R,\mm)$ be a $2$-dimensional local \CM ring and let $x,y\in\mm$ be any system of parameters of $R.$ Let $I$ be any ideal integral over $(x,y).$ Then for all $n\geq 1,$
\beqn
\lm\lf\frac{I^{n+1}}{(x,y)I^n}\rg-\lm\lf\frac{(I^n:(x,y))}{I^{n-1}}\rg = \Delta^2\lf{P_I(n-1)-H_I(n-1)}\rg.
\eeqn 
\el S. Huckaba \cite{huckaba} extended this result for dimension $d\geq 1.$ In section three, for $l\geq 1,$ $1\leq k\leq d$ and  $\underline{y}^{[l]} = ({y_1}^l,\ldots,{y_k}^l),$ we introduce an analogue of the  Kirby-Mehran complex \cite{KM} for multigraded filtrations: 
$\lmc :$ $$0\longrightarrow {\frac{R}{{\fil}(\n)}}\overset{d_k}\longrightarrow\lf{\frac{R}{{\fil}(\n+le)}}\rg^{k \choose {1}}\overset{d_{k-1}}\longrightarrow\cdots\overset{d_2}\longrightarrow\lf{\frac{R}{{\fil}(\n+(k-1)le)}}\rg^{{k}\choose {k-1}}\overset{d_1}\longrightarrow {\frac{R}{{\fil}(\n+kle)}}\overset{d_0}\longrightarrow 0$$ 
and prove an analogue of Huneke's Fundamental Lemma for multigraded filtration of ideals.
\bt
Let $(R,\mm)$ be a \CM local ring of dimension $d\geq 1$ with infinite residue field, $\idl$ be $\mm$-primary ideals of $R$ and $\fil=\lbrace\fil(\n)\rbrace_{\n\in \ZZ^s}$ be an $\I$-admissible filtration of ideals in $R.$ Let $\mathcal{A}={\lbrace}{x_{ij}\in{I_i}: j=1,\ldots,d;\ii}\rbrace$ be any complete reduction of $\fil,$ $y_j=x_{1j}\cdots x_{sj}$ for all $j=1,\ldots,d.$ Let $\underline y=\yy$ and $J=(\underline y).$ Then for all $\n\in\ZZ^s,$ $${\Delta}^d(\po(\n)-\ho(\n))=\lm\lf\frac{{\fil}(\n+de)}{J{\fil}(\n+(d-1)e)}\rg-\sum_{i=2}^{d}{(-1)^i}{\lm({H_i(C.(\underline y,\fil(\n))}))}.$$\et
In section four, for $d\geq 2,$ we compute $\n$ degree component of local cohomolgy module $H_{\R_{++}}^1(\R(\fil))_{\n}$ for all $\n\in\NN^s$ and give an equivalent criterion for vanishing of $H_{\R_{++}}^1(\R(\fil))_{\n}$ for all $\n\in\NN^s.$ We discuss vanishing of Hilbert coefficients and generalize some results due to T. Marley \cite{marleythesis} in \CM local ring of dimension $1\leq d\leq 2$.

\bt
Let $(R,\mm)$ be a \CM local ring of dimension $1\leq d\leq 2$ with infinite residue field, $\idl$ be $\mm$-primary ideals of $R$ and $\fil=\lbrace\fil(\n)\rbrace_{\n\in \ZZ^s}$ be an $\I$-admissible filtration of ideals in $R.$ Let $e_{(d-1)e_i}(\fil)=0$ for $\ii.$ Then 
\ben
\item For $d=1,$ $\po(\n)=\ho(\n)$ for all $\n\in\NN^s.$
\item For $d=2,$ if $H_{\R_{++}}^1(\R(\fil))_{\n}=0$ for all $\n\in\NN^s$ then $\po(\n)=\ho(\n)$ for all $\n\in\NN^s$ and $e_{\underline{0}}(\fil)=0.$\een\et
\bt
Let $(R,\mm)$ be a \CM local ring of dimension $1\leq d\leq 2$ with infinite residue field and $I,J$ be $\mm$-primary ideals of $R.$ Let $\fil=\lbrace\fil(r,s)\rbrace_{r,s\in \ZZ}$ be a $\ZZ^2$-graded $(I,J)$-admissible filtration of ideals in $R.$ Then the following statements are equivalent.\ben
\item $e_{(d-1)e_i}(\fil)=0$ for $i=1,2.$
\item $I$ and $J$ are generated by system of parameters, $\po(r,s)=\ho(r,s)$ for all $r,s\in\NN$ and $\fil(r,s)=I^rJ^s$ for all $r,s\in\ZZ.$ 
\item $e_{\alpha}(\fil)=0$ for $|\alpha|\leq d-1.$
\een\et
In \cite{marley}, Marley proved that for \CM local ring of dimension $d\geq 1,$ $r(I)=n(I)+d$ under some depth condition of the associated graded ring of $I.$ In his thesis \cite{marleythesis}, Marley extended this result for $\ZZ$-graded $I$-admissible filtrations. We generalize this result for multigraded filtration of ideals when $d=1,2.$ In section five, we prove the following theorem.
\bt
Let $(R,\mm)$ be a \CM local ring of dimension one with infinite residue field and $\idl$ be $\mm$-primary ideals of $R.$ Let $\fil=\lbrace\fil(\n)\rbrace_{\n\in \ZZ^s}$ be an $\I$-admissible filtration of ideals in $R$ and $\mathcal {A}=\lbrace{a_i\in{I_i}:\ii}\rbrace$ be a complete reduction of $\fil.$ Then $$\PP(\fil)\subseteq\NN^s\mbox{  and  } \PP(\fil)={\mathcal{R}}_{\mathcal A}(\fil).$$ Moreover, the set $\R_{\mathcal A}(\fil)$ is independent of any complete reduction $\mathcal A$ of $\fil.$
\et
We also show that for one-dimensional \CM local ring $(R,\mm),$ $r_{\mathcal A}(\fil)$ is independent of any complete reduction $\mathcal A$ of $\fil.$
\\In section six, we provide a relation between reduction vectors of good complete reductions and postulation vectors of multigraded filtration of ideals in two-dimensional \CM local rings. For bigraded filtration we prove a result which relates the Cohen-Macaulayness of the bigraded Rees algebra, the complete reduction number, reduction numbers and the joint reduction number.
\bt
 Let $(R,\mm)$ be a \CM local ring of dimension two with infinite residue field and $\idl$ be $\mm$-primary ideals of $R$ and $s\geq 2.$ Let $\fil=\lbrace\fil(\n)\rbrace_{\n\in \ZZ^s}$ be an $\I$-admissible filtration of ideals in $R$ and $\mathcal {A}={\lbrace}{x_{ij}\in{I_i}: j=1,2;\ii}\rbrace$ be a \gcrr of $\fil.$ Let $H_{\R_{++}}^1(\R(\fil))_{\n}=0$ for all $\n\geq\underline 0.$ Then $\PP(\fil)\subseteq\NN^s$ and
there exists a one-to-one correspondence 
$$f:{\displaystyle{\Large{\Large{\PP(\fil)\longleftrightarrow\lbrace{\rrr\in{\R_{\mathcal A}(\fil)}\mid\rrr\geq e}\rbrace}}}}$$ defined by $f(\n)=\n+e$ where $f^{-1}(\rrr)=\rrr-e.$ 
\et
\bt
Let $(R,\mm)$ be a \CM local ring of dimension two with infinite residue field and $\idl$ be $\mm$-primary ideals of $R$ and $s\geq 2.$ Let $\fil=\lbrace\fil(\n)\rbrace_{\n\in \ZZ^s}$ be an $\I$-admissible filtration of ideals in $R$ and $H_{\R_{++}}^1(\R(\fil))_{\n}=0$ for all $\n\geq\underline 0.$ Then the following statements are equivalent.
\ben
{
\item[(1)] $\PP(\fil)=\NN^s,$ i.e. $\po(\n)=\ho(\n)$ for all $\n\geq\underline 0.$ 
\item[(2)] $r_{\mathcal A}(\fil)\leq {1}$ for any \gcrr $\mathcal {A}$ of $\fil.$
\item[$(2')$] There exists a \gcrr $\mathcal {A}$ of $\fil$ such that $r_{\mathcal A}(\fil)\leq {1}.$
}
\een \et
In order to state the final result, we recall the definition of joint reduction of multigraded filtrations \cite{msv}. The {\bf{joint reduction}} of $\fil$ of type ${\bf{q}}=(q_1,\ldots,q_s)\in \NN^s$ is a collection of $q_i$ elements $x_{i1},\ldots,x_{iq_i}\in I_i$ for all $\ii$ such that $q_1+\cdots+q_s=d$ and $$\sum_{i=1}^{s}{\sum_{j=1}^{q_i}{x_{ij}\fil(\n-e_i)}}=\fil(\n)\mbox{ for all large }\n.$$ We say the {\bf{joint reduction number of $\fil$ with respect to a joint reduction $\lbrace x_{ij}\in I_i:j=1,\ldots,q_i;\ii\rbrace$ of type $q$ is zero}} if 
$$\sum_{i=1}^{s}{\sum_{j=1}^{q_i}{x_{ij}\fil(\n-e_i)}}=\fil(\n) \mbox{ for all } 
\n \geq \sum_{i\in A} e_i, \mbox{ where } A=\{i | q_i \neq 0\}.$$
We say the {\bf{joint reduction number of $\fil$ of type $q$ is zero}} if the joint reduction number of $\fil$ with respect to any joint reduction of type $q$ is zero.
\begin{theorem}
Let $(R,\mm)$ be a \CM local ring of dimension two with infinite residue field and $I,J$ be $\mm$-primary ideals of $R.$ Let $\fil=\lbrace\fil(\n)\rbrace_{\n\in \ZZ^2}$ be a $\ZZ^2$-graded $(I,J)$-admissible filtration of ideals in $R.$ Then the following are equivalent.
\ben
{\item[(1)] The Rees algebra $\R(\fil)$ is Cohen-Macaulay.
\item[(2)] $\PP(\fil)=\NN^2,$ i.e. $\po(\n)=\ho(\n)$ for all $\n\geq\underline 0.$  
\item[(3)] For any \gcrr $\mathcal {A}$ of $\fil,$ $r_{\mathcal A}(\fil)\leq {1}$ and $H_{\R_{++}}^1(\R(\fil))_{\n}=0$ for all $\n\geq\underline 0.$
\item[$(3')$] There exists a \gcrr $\mathcal {A}$ of $\fil$ such that $r_{\mathcal A}(\fil)\leq {1}$ and $H_{\R_{++}}^1(\R(\fil))_{\n}=0$ for all $\n\geq\underline 0.$
\item[(4)] For the filtrations $\mathcal F^{(i)}=\lbrace\fil(ne_i)\rbrace_{n\in \ZZ},$ $r(\mathcal F^{(i)}) \leq 1$ where $i=1,2$ and the joint reduction number of $\mathcal F$ of type $e$ is zero. 
}
\een\end{theorem}
 \section{Preliminary results}
 In this section we discuss the existence of good complete reduction of an $\I$-admissible filtration $\fil=\lbrace\fil(\n)\rbrace_{\n\in \ZZ^s}$ and prove some results about Hilbert coefficients which we need in the following sections. For an admissible multigraded filtration  $\fil=\lbrace\fil(\n)\rbrace_{\n\in{\ZZ}^s}$, by Rees' Lemma \cite[Lemma 1.2]{rees3}, \cite[Lemma 2.2]{msv}, we get elements $x_i\in I_i$ for all $\ii,$ called superficial elements for $\fil$, such that for each $i,$ there exists an integer $r_i$ and $(x_i)\cap\fil(\n)=x_i\fil(\n-e_i)$ for all $\n\geq r_ie_i.$ In \cite[Theorem 1.3]{rees3}, Rees proved existence of complete reduction of the filtration $\lbrace\I^{\n}\rbrace_{\n\in \ZZ^s}.$ Using the same lines of proof of this theorem and existence of superficial elements we get the following theorem.
\bt
Let $(R,\mm)$ be a Cohen-Macaulay local ring of dimension $d$ with infinite residue field and $\idl$ be $\mm$-primary ideals of $R.$ Let $\fil$ be an $\I$-admissible filtration of ideals in $R.$ Then there exists a \gcrr of $\fil.$
\et
\bl\label{one}
Let $(R,\mm)$ be a Noetherian local ring of dimension $d\geq 1$ and $\idl$ be $\mm$-primary ideals of $R.$ Put $I=\id.$ Then $$
\displaystyle\sum\limits_{\substack{\alpha=({\alpha}_1,\ldots,{\alpha}_s)\in{\NN}^s \\ |\alpha|= d}}\frac{{d!}e_{\alpha}(\I)}{{\alpha}_1!\cdots{{\alpha}_s}!}={e_0(I)}\mbox{ and }e_{\underline{0}}(\I)=e_d(I).$$
\bpf Since for large $n,$ 
$$\pf(ne)=\hf(ne)=\lm\lf\frac{R}{{\I}^{ne}}\rg=\lm\lf\frac{R}{I^n}\rg=P_I{(n)},$$ 
comparing the coefficients of $n^d$ and constant terms we get the required result.
\eepf
\el
\bp\label{math1}
Let $s\geq 1$ be a fixed integer and $i_1,\ldots,i_s\in\NN$ such that $g=i_1+\cdots+i_s\geq 1.$ Then ${\Delta}^g({n_1}^{i_1}\cdots{n_s}^{i_s})=g!$ where $n_k\in\ZZ$ for all $k=1,\ldots,s.$
\bpf
We use induction on $g.$ Let $g=1.$ Then \wlg assume $i_1=1$ and $i_k=0$ for all $k\neq 1.$ Therefore 
${\Delta}^1(n_1)=(n_1+1)-n_1=1.$ Let $g\geq 2$ and assume the result is true up to $g-1.$ Now
\beqn
{\Delta}^g({n_1}^{i_1}\cdots{n_s}^{i_s})&=&{\Delta}^{g-1}\Big[{\Delta}^1({n_1}^{i_1}\cdots{n_s}^{i_s})\Big]\\&=& {\Delta}^{g-1}\Big[(n_1+1)^{i_1}\cdots(n_s+1)^{i_s}- {n_1}^{i_1}\cdots{n_s}^{i_s}\Big]\\&=& {\Delta}^{g-1}\Big[\sum_{k=1}^{s}{i_k}({n_1}^{i_1}\cdots{n_k}^{i_k-1}\cdots{n_s}^{i_s})\Big]\\&=&\sum_{k=1}^{s}{i_k}{\Delta}^{g-1}({n_1}^{i_1}\cdots{n_k}^{i_k-1}\cdots{n_s}^{i_s})\\&=&\sum_{k=1}^{s}{i_k}(g-1)!\\&=& g!.
\eeqn
\eepf
\ep
\bp\label{kolkata}
Let $(R,\mm)$ be a Noetherian local ring of dimension $d\geq 1$ and $\idl$ be $\mm$-primary ideals of $R.$ Let $\fil=\lbrace\fil(\n)\rbrace_{\n\in \ZZ^s}$ be an $\I$-admissible filtration of ideals in $R.$ Then 
\ben
{
\item[(1)] $e_{\alpha}(\fil)=e_{\alpha}(\I)$ for all $\alpha\in\NN^s$ where $|\alpha|=d.$
\item[(2)]${\Delta}^d(\po(\n))={\Delta}^d(\pf(\n))=e_0(\id).$
\item[(3)] For an ideal $J\subseteq\fil(e)$ such that $J{\fil(\n)}=\fil(\n+e)$ for all large $\n\in\NN^s,$ we have ${\Delta}^d(\po(\n))=e_0(J)=e_0(\id).$
}
\een
\bpf
$(1)$ This follows from \cite[Theorem 2.4]{rees3}.
\\$(2)$ Using Proposition \ref{math1}, we get
\beqn
{\Delta}^d(\po(\n))&=&{\Delta}^d\lf\displaystyle\sum_{{|\alpha|}\leq d}(-1)^{d-|\alpha|}{e_{\alpha}(\fil)}\binom{{n_1}+{{\alpha}_1}-1}{{\alpha}_1}\cdots\binom{{n_s}+{{\alpha}_s}-1}{{\alpha}_s}\rg\\&=&{\Delta}^d\lf\displaystyle\sum_{|\alpha|=d}\frac{e_{\alpha}(\fil)}{{{{\alpha}_1}!}\cdots{{{\alpha}_s}!}}{n_1}^{{\alpha}_1}\cdots{n_s}^{{\alpha}_s}\rg\\&=&\displaystyle\sum_{{|\alpha|}=d}\frac{e_{\alpha}(\fil)}{{{{\alpha}_1}!}\cdots{{{\alpha}_s}!}}\Big[{\Delta}^d({n_1}^{{\alpha}_1}\cdots{n_s}^{{\alpha}_s})\Big]\\&=&\displaystyle\sum_{{|\alpha|}=d}\frac{e_{\alpha}(\fil)}{{{{\alpha}_1}!}\cdots{{{\alpha}_s}!}}({\alpha}_1+\cdots+{\alpha}_s)!\\&=&\displaystyle\sum_{|\alpha|=d}d!\frac{e_{\alpha}(\fil)}{{{{\alpha}_1}!}\cdots{{{\alpha}_s}!}}.
\eeqn In a similar way we get ${\Delta}^d(\pf(\n))= \displaystyle\sum_{|\alpha|=d}d!\frac{e_{\alpha}(\I)}{{{{\alpha}_1}!}\cdots{{{\alpha}_s}!}}.$ By part $(1)$ and Lemma \ref{one}, ${\Delta}^d(\po(\n))={\Delta}^d(\pf(\n))=e_0(\id).$
\\$(3)$ Since $J{\fil(\n)}=\fil(\n+e)$ for all large $\n,$ there exists an integer $k\in\NN$ such that $J{\fil(\n)}=\fil(\n+e)$ for all $\n\geq ke.$ Now for all $n\geq {k},$ $J^{n-k}{\fil(ke)}=\fil(ne)$ and hence $J^n\subseteq\fil(ne)\subseteq J^{n-k}.$ This implies $$\lm\lf\frac{R}{{J}^{n-k}}\rg\leq\lm\lf\frac{R}{{\fil}(ne)}\rg\leq\lm\lf\frac{R}{{J}^{n}}\rg$$ for all $n\geq k.$ Therefore for all $n\geq k,$ we have $$\lim_{n\to\infty}\frac{P_J(n-k)}{n^d/d!}\leq\lim_{n\to\infty}\frac{\po(ne)}{n^d/d!}\leq \lim_{n\to\infty}\frac{P_J(n)}{n^d/d!}$$ which implies ${\Delta}^d(\po(\n))=e_0(J).$ Hence using part $(2),$ we get the required result.

\eepf
\ep
\bp\label{NH1}
Let $(R,\mm)$ be a Noetherian local ring of dimension $d\geq 1,$ $\depth R\geq 1$ and $\idl$ be $\mm$-primary ideals of $R.$ Let $\fil=\lbrace\fil(\n)\rbrace_{\n\in \ZZ^s}$ be an $\I$-admissible filtration of ideals in $R$ and $\mathcal {A}={\lbrace}{x_{ij}\in{I_i}: j=1,\ldots,d;\ii}\rbrace$ be a \crr of $\fil.$ Then $J_i=(x_{i1},\ldots,x_{id})$ is a reduction of $I_i$ for all $\ii.$
\bpf
For all large $\n,$ $J{\fil(\n)}=\fil(\n+e)$ where $y_j=x_{1j}\cdots x_{sj}$ for all $j=1,\ldots,d$ and $J=(\yy).$ Since $\fil$ is an $\I$-admissible filtration, for each $i\in\{1,\dots,s\},$ there exists $r_i\in\NN$ such that for all $\n\geq r_ie_i,$ $\fil(\n+e_i)=I_i\fil(\n).$ Hence for all large $\n,$
$$
I_i{\fil}(\n+e-e_i)\supseteq J_i{\fil}(\n+e-e_i)\supseteq  J{\fil}(\n)={\fil}(\n+e)\supseteq  I_i{\fil}(\n+e-e_i).
$$ Now by \cite[Lemma 1.5]{rees3}, $J_i$ is reduction of $I_i$ for all $\ii.$
\eepf
\ep

\bp\label{a}
Let $f(\n):{\ZZ}^s\rightarrow {\ZZ}$ be an integer valued function such that for all large $\n\in\NN^s,$ $f(\n)=0.$ Let  $\mathcal{B}=\lbrace{\n\in{\ZZ}^s:f(\m)=0\mbox{ for all } \m\geq\n}\rbrace$ and $\mathcal{C}_j=\lbrace{\n\in{\ZZ}^s:{\Delta}^{j}(f(\m))=0\mbox{ for all }\m\geq\n }\rbrace$ for all $j\geq 0.$ Then for $j\geq 0,$ $\mathcal{B}=\mathcal{C}_j.$
\bpf
For $j=0,$ the result is true due to the definition of $\mathcal{B}.$ It is enough to prove the statement for $j=1.$ Let $\n\in\mathcal{B}.$ Then $\Delta^{1}(f(\m))=f(\m+e)-f(\m)=0$ for all $\m\geq\n.$ This implies $\n\in{\mathcal{C}_1}.$ Conversely, let $\n\in{\mathcal{C}_1}.$ Then for all $\m\geq\n,$ $0=\Delta^{1}(f(\m))=f(\m+e)-f(\m).$ Let $\kr\in{\NN}^s$ such that $f(\rrr)=0$ for all $\rrr\geq\kr.$ Let $\m\geq\n.$ For all $\ii,$ define 
\[ u(m)_i = \left\{
  \begin{array}{l l}
   k_i-m_i  & \quad \text{if $k_i>m_i$ }\\
    0 & \quad \text{if $k_i\leq m_i$ }
  \end{array} \right.\] Let $u(m)=\max\lbrace{u(m)_1,\ldots,u(m)_s}\rbrace+1.$ Then for all $\m\geq\n,$ $$0=f(\m+u(m)e)=\cdots=f(\m+e)=f(\m).$$ Hence $\n\in\mathcal{B}.$
\eepf\ep
\bp\label{r}
Let $R=\bigoplus\limits_{\n\in{\NN}^s}R_{\n}$ be an ${\NN}^s$-graded ring, $S$ be an ${\NN}^s$-graded $R$-algebra and $b\in R_e.$ Let $S_{\n}\overset{.b}\longrightarrow S_{\n+e}$ be an injective map for all large $\n$ and grade$(S_{++})\geq 1.$ Then $b$ is a nonzerodivisor of $S.$
\bpf
Let $\m\in{\NN}^s$ such that for all $\n\geq\m,$ $S_{\n}\overset{.b}\longrightarrow S_{\n+e}$ is an injective map. Let $x\in {(0:_{S}b)\cap S_{\kr}}$ for some $\kr\in\NN^s.$ We show that $x(S_{++})^{m+1}=0$ where $m=\max\lbrace{m_1,\ldots,m_s}\rbrace.$ Let $0\neq z\in {(S_{++})^{m+1}\cap S_{\underline p}}.$ Now $xz\in S_{\kr+\underline p}$ and $bxz=0.$ Since $\kr+\underline p\geq (m+1)e,$ $xz=0.$ Thus $x\in (0:_{S}{(S_{++})^{m+1}})=0.$
\eepf
\ep
\section{An analogue of Huneke's fundamental lemma}
Throughout this section $(R,\mm)$ is a \CM local ring of dimension $d\geq 1$ with infinite residue field and $\idl$ are $\mm$-primary ideals of $R.$ Let $\fil=\lbrace\fil(\n)\rbrace_{\n\in \ZZ^s}$ be an $\I$-admissible filtration of ideals in $R.$ In \cite{jv}, A. V. Jayanthan and J. K. Verma generalized the Kirby-Mehran complex \cite{KM}, \cite{huckaba-marley} for the bigraded filtration $\{I^rJ^s\}_{r,s\in\ZZ}$ where $I,J$ are $\mm$-primary ideals and studied the relation between cohomology modules of this complex and local cohomology modules of $\R(I,J).$ We construct a multigraded analogue of the Kirby-Mehran complex and compute its homology modules. As a consequence of this we prove an analogue of Huneke's Fundamental Lemma \cite{huneke},\cite{huckaba}.  \\Let $y_1,\ldots,y_k$ be elements in $\id$ where $1\leq k\leq d.$ For $l\geq 1$ and $({\underline {yt}})^{[l]}={y_1}^l\underline{t}^{le},\ldots,{y_k}^l\underline{t}^{le},$ consider the Koszul complex $K.((\underline{yt})^{[l]},\mathcal{R}'(\fil)):$ $$0\longrightarrow\mathcal{R}'(\fil)\longrightarrow\mathcal{R}'(\fil)(le)^{k \choose {1}}\longrightarrow\cdots\longrightarrow\mathcal{R}'(\fil)((k-1)le)^{k \choose {k-1}}\longrightarrow \mathcal{R}'(\fil)(kle)\longrightarrow 0.$$ This complex has a ${\ZZ}^s$-graded structure inherited from $\mathcal{R}'(\fil).$ The graded component of degree $\n$ of the above Koszul complex is $K._{\n}((\underline{yt})^{[l]},\mathcal{R}'(\fil)):$ $$0\longrightarrow {{\fil}(\n)}\longrightarrow(\fil(\n+le))^{k \choose {1}}\longrightarrow\cdots\longrightarrow({{{\fil}(\n+(k-1)le)}})^{{k}\choose {k-1}}\longrightarrow {{{\fil}(\n+kle)}}\longrightarrow 0.$$ Let $\underline{y}^{[l]} = {y_1}^l,\ldots,{y_k}^l.$ We can consider the above complex as a subcomplex of the Koszul complex $K.({y}^{[l]},R):$ $$0\longrightarrow R\longrightarrow R^{{k}\choose {1}}\longrightarrow\cdots\longrightarrow R^{{k}\choose {k-1}}\longrightarrow R\longrightarrow 0.$$ Hence we have a chain map of complexes $0\longrightarrow K._{\n}((\underline{yt})^{[l]},\mathcal{R}'(\fil))\longrightarrow K.(\underline{y}^{[l]},R)$ which produces a quotient complex $\lmc :$ $$0\longrightarrow {\frac{R}{{\fil}(\n)}}\overset{d_k}\longrightarrow\lf{\frac{R}{{\fil}(\n+le)}}\rg^{k \choose {1}}\overset{d_{k-1}}\longrightarrow\cdots\overset{d_2}\longrightarrow\lf{\frac{R}{{\fil}(\n+(k-1)le)}}\rg^{{k}\choose {k-1}}\overset{d_1}\longrightarrow {\frac{R}{{\fil}(\n+kle)}}\overset{d_0}\longrightarrow 0.$$ 
In the following proposition we compute homology modules of the above complex. 
\bp\label{homology}
For all $l\geq 1,\mbox{ } \n\in{\ZZ}^s\mbox{ and } 1\leq k\leq d,$
\\ $(1)$ $H_0(\lmc)=\displaystyle{\frac{R}{({\fil}(\n+kle),\underline{y}^{[l]})}},$
\\ $(2)$ $H_k(\lmc)=\displaystyle{\frac{({\fil}(\n+le):(\underline{y}^{[l]}))}{{\fil}(\n)}},$
\\ $(3)$ if $\yk$ is a regular sequence then $H_1(\lmc)=\displaystyle{\frac{({(\underline{y}^{[l]})}\cap{{\fil}(\n+kle)})}{(\underline{y}^{[l]}){\fil}(\n+(k-1)le)}}.$
\ep
\bpf
$(1)$ Since $\ker d_0={R/{\fil}(\n+kle)}$ and $\im d_1=\lf(\underline{y}^{[l]})+\fil(\n+kle)\rg /{\fil(\n+kle)},$ we get 
$$H_0(\lmc)=\frac{\ker d_0}{\im d_1}=\frac{R}{({\fil}(\n+kle),\underline{y}^{[l]})}.$$
$(2)$ Since $\im d_{k+1}=0,$\beqn
H_k(\lmc)&=& \ker\hspace{0.1cm} d_k\\&=&\lbrace{x+{\fil}(\n)\in {R/{\fil}(\n)}|x{y_i}^l\in{{\fil}(\n+le)}\mbox{ for all }i=1,\ldots,k}\rbrace\\&=&\lbrace{x+{\fil}(\n)\in {R/{\fil}(\n)}|x\in{\cap}_{i=1}^{k}({{\fil}(\n+le)}:({y_i}^l))}\rbrace\\&=&\lbrace{x+{\fil}(\n)\in {R/{\fil}(\n)}|x\in{({\fil}(\n+le)}:(\underline{y}^{[l]}))}\rbrace.
\eeqn
$(3)$ Since $\yk$ is a regular sequence, the following sequence is exact $$R^{{k}\choose{k-2}}\overset{\phi_2}\longrightarrow R^{{k}\choose{k-1}}\overset{\phi_1}\longrightarrow (\underline{y}^{[l]})\overset{\phi_0}\longrightarrow 0.$$ Tensoring by $\frac{R}{\fil(\n+(k-1)le)},$ we get an exact sequence $$\lf{\frac{R}{\fil(\n+(k-1)le)}}\rg^{{k}\choose{k-2}}\overset{\overline{\phi}_2}\longrightarrow \lf{\frac{R}{\fil(\n+(k-1)le)}}\rg^{{k}\choose{k-1}}\overset{\overline{\phi}_1}\longrightarrow \frac{(\underline{y}^{[l]})}{(\underline{y}^{[l]})\fil(\n+(k-1)le)}\overset{\overline{\phi}_0}\longrightarrow 0.$$ Hence $\im {\overline{\phi}_2}=\im {d_2}$ and we get the following commutative diagram of exact rows:
$$\displaystyle{
\begin{CD}
0 @>>>{\im {d_2}}@>>>{\lf{\frac{R}{\fil(\n+(k-1)le)}}\rg^{{k}\choose{k-1}}}@>>>{\frac{(\underline{y}^{[l]})}{(\underline{y}^{[l]})\fil(\n+(k-1)le)}}@>>> 0\\
@. @ViVV @VidVV @V{\theta}VV @.\\
0 @>>>{\ker {d_1}}@>>>{\lf{\frac{R}{\fil(\n+(k-1)le)}}\rg^{{k}\choose{k-1}}} @>>>{\frac{R}{\fil(\n+kle)}}
\end{CD}}$$ where $i$ is the inclusion map and $id$ is the identity map. Then by the Snake Lemma, $$H_1(\lmc)=\frac{\ker {d_1}}{\im {d_2}}\cong \ker\theta = \frac{({(\underline{y}^{[l]})}\cap{{\fil}(\n+kle)})}{(\underline{y}^{[l]}){\fil}(\n+(k-1)le)}.$$
\eepf

\bt\label{use}{(Analogue of Huneke's Fundamental Lemma) }
Let $(R,\mm)$ be a \CM local ring of dimension $d\geq 1$ with infinite residue field, $\idl$ be $\mm$-primary ideals of $R$ and $\fil=\lbrace\fil(\n)\rbrace_{\n\in \ZZ^s}$ be an $\I$-admissible filtration of ideals in $R.$ Let $\mathcal{A}={\lbrace}{x_{ij}\in{I_i}: j=1,\ldots,d;\ii}\rbrace$ be any complete reduction of $\fil,$ $y_j=x_{1j}\cdots x_{sj}$ for all $j=1,\ldots,d.$ Let $\underline y=\yy$ and $J=(\underline y).$ Then for all $\n\in\ZZ^s,$ $${\Delta}^d(\po(\n)-\ho(\n))=\lm\lf\frac{{\fil}(\n+de)}{J{\fil}(\n+(d-1)e)}\rg-\sum_{i=2}^{d}{(-1)^i}{\lm({H_i(C.(\underline y,\fil(\n))}))}.$$
\bpf  By Propositions \ref{kolkata} and \ref{homology}, we get
\beqn
{\Delta}^d(\po(\n)-\ho(\n))&=& e_0(J)-{\Delta}^d(\ho(\n))\\&=& e_0(J)-\sum_{i=0}^{d}{(-1)^i}{\lm({H_i(C.(\underline y,\fil(\n))}))}\\&=&\lm{\lf{\frac{R}{(\underline y)}}\rg}-\lm{\lf{\frac{R}{({\fil}(\n+de),{(\underline y)})}}\rg}+\lm\lf\frac{{(\underline y)}\cap{{\fil}(\n+de)}}{(\underline y){\fil}(\n+(d-1)e)}\rg\\&&-\sum_{i=2}^{d}{(-1)^i}{\lm({H_i(C.(\underline y,\fil(\n))}))}\\&=&\lm\lf\frac{{\fil}(\n+de)}{(\underline y){\fil}(\n+(d-1)e)}\rg-\sum_{i=2}^{d}{(-1)^i}{\lm({H_i(C.(\underline y,\fil(\n))}))}\\&=&\lm\lf\frac{{\fil}(\n+de)}{J{\fil}(\n+(d-1)e)}\rg-\sum_{i=2}^{d}{(-1)^i}{\lm({H_i(C.(\underline y,\fil(\n))}))}.
\eeqn
\eepf
\et
\section{Vanishing of Hilbert coefficients}
In this section we compute the local cohomology module $H_{\R_{++}}^1(\R(\fil))_{\n}$ for all $\n\geq\underline 0$ \cite{blancafort},\cite{jv}. We discuss the vanishing of Hilbert coefficients of an $\I$-admissible filtration $\fil$ and generalize some results due to Marley \cite{marleythesis} in \CM local ring of dimension $1\leq d\leq 2.$

For a $\ZZ^s$-graded $\I$-filtration $\fil=\lbrace\fil(\n)\rbrace_{\n\in \ZZ^s},$ the filtration $\breve{\fil}=\{\breve \fil(\n)\}_{\n\in\ZZ^s}$ of ideals is called Ratliff-Rush closure filtration of $\fil$ where $\breve \fil(\n)=\bigcup\limits_{k\geq 1} (\fil(\n+ke):\fil(e)^k)$ for all $\n\in\NN^s$ and $\breve \fil(\n)=\breve\fil(\n^+)$ for all $\n\in\ZZ^s$ \cite{msv}. To compute $H_{\R_{++}}^1(\R(\fil))_{\n}$ we follow the lines of the proof \cite[Theorem 3.3]{blancafort}. 
\bp\label{h1}
Let $(R,\mm)$ be a \CM local ring of dimension $d\geq 2$ with infinite residue field, $\idl$ be $\mm$-primary ideals of $R$ and $\fil=\lbrace\fil(\n)\rbrace_{\n\in \ZZ^s}$ be an $\I$-admissible filtration of ideals in $R.$ Then for all $\n\in\NN^s,$ $$H_{\R_{++}}^1(\R(\fil))_{\n}\cong\frac{\breve\fil(\n)}{\fil(\n)}.$$
\bpf
Let $\mathcal{A}={\lbrace}{x_{ij}\in{I_i}: j=1,\ldots,d;\ii}\rbrace$ be any complete reduction of $\fil,$ $y_j=x_{1j}\cdots x_{sj}$ for all $j=1,\ldots,d.$ For each $\n,$ $l\geq 1,$ $({\underline {yt}})^{[l]}={y_1}^l\underline{t}^{le},\ldots,{y_d}^l\underline{t}^{le}$ and ${\underline y}^{[l]}={y_1}^l,\ldots,{y_d}^l,$ we have the following exact sequence 
$$0\longrightarrow K._{\n}(({\underline {yt}})^{[l]},\mathcal{R}'(\fil))\longrightarrow K.({\underline y}^{[l]},R)\longrightarrow C.({\underline y}^{[l]},\fil(\n)) \longrightarrow 0.$$  
For each $i\in\lbrace 1,\ldots,d\rbrace$ the following commutative diagram of complexes :
$$\begin{CD} K.({y_i}^l,R):~~~~0@>>>R@>{{y_i}^l}>>R@>>>0\\
@. @VidVV @V{y_i}VV  @.\\
K.({y_i}^{l+1},R):~~0@>>>R@>{{y_i}^{l+1}}>>R@>>>0 \end{CD}$$ gives a map $K.({y_1}^l,\ldots,{y_d}^l,R)=\otimes_{i=1}^{d} K.({y_i}^l,R)\rightarrow K.({y_1}^{l+1},\ldots,{y_d}^{l+1},R)=\otimes_{i=1}^{d} K.({y_i}^{l+1},R).$
The maps can be restricted to $K._{\n}({y_1}^l\underline{t}^{le},\ldots,{y_d}^l\underline{t}^{le},\mathcal{R}'(\fil)).$ Hence for all $l\geq 1,$ we get morphisms of exact sequences 
$$
\begin{CD} 0@>>>K._{\n}(({\underline {yt}})^{[l]},\mathcal{R}'(\fil))@>>>K.({\underline y}^{[l]},R)@>>>C.({\underline y}^{[l]},\fil(\n))@>>>0\\
@. @VVV @VVV @VVV  @.\\
0@>>>K._{\n}(({\underline {yt}})^{[l+1]},\mathcal{R}'(\fil))@>>>K.({\underline y}^{[l+1]},R)@>>>C.({\underline y}^{[l+1]},\fil(\n))@>>>0
 \end{CD}$$ which produce an inductive system of exact sequences of complexes.
Applying $\displaystyle\lim_{\stackrel{\longrightarrow}{l}}$ to the long exact sequence of cohomology modules, we get  
 $$0\longrightarrow H_{({\underline {yt}})}^0(\mathcal{R}'(\fil))_{\n}\longrightarrow H_{(\underline y)}^0(R)\longrightarrow\displaystyle\lim_{\stackrel{\longrightarrow}{l}} H^0(C.({\underline y}^{[l]},\fil(\n)))\longrightarrow H_{({\underline {yt}})}^1(\mathcal{R}'(\fil))_{\n}\longrightarrow\cdots.$$
Since $(R,\mm)$ is Cohen-Macaulay, $ H_{(\underline y)}^i(R)=0$ for $0\leq i\leq d-1.$ Hence $$H_{({\underline {yt}})}^1(\mathcal{R}'(\fil))_{\n}\cong \displaystyle\lim_{\stackrel{\longrightarrow}{l}} H^0(C.({\underline y}^{[l]},\fil(\n)))=\displaystyle\lim_{\stackrel{\longrightarrow}{l}} H_d(C.({\underline y}^{[l]},\fil(\n)))=\displaystyle\lim_{\stackrel{\longrightarrow}{l}}\frac{({\fil}(\n+le):(\underline{y}^{[l]}))}{{\fil}(\n)}.$$ Since $\sqrt{\R_{++}}=\sqrt{({\underline {yt}})},$ by \cite[Proposition 3.1, Proposition 4.2]{msv}, for all $\n\in\NN^s,$  $$H_{\R_{++}}^1(\R(\fil))_{\n}\cong H_{\R_{++}}^1(\R'(\fil))_{\n}\cong\frac{\breve\fil(\n)}{\fil(\n)}.$$
\eepf
\ep
For all $\ii,$ we denote the {\bf{associated multigraded ring of $\fil$ \wrt $\fil(e_i)$}} by $\displaystyle\ggi=\bigoplus\limits_{\n\in{\NN}^s}\frac{{\fil}(\n)}{{\fil}(\n+e_i)}.$ For $\fil=\{\I^{\n}\}_{\n \in \ZZ^s},$ we set $\ggi=G_{i}(\I).$ In the next proposition we give an equivalent criterion for the vanishing of $H_{\R_{++}}^1(\R(\fil))_{\n}$ for all $\n\geq \underline 0$ in terms of grade$(\gi)$ for all $\ii.$ 
\bp\label{grade}
Let $(R,\mm)$ be a \CM local ring of dimension $d\geq 2$ with infinite residue field, $\idl$ be $\mm$-primary ideals of $R$ and $\fil=\lbrace\fil(\n)\rbrace_{\n\in \ZZ^s}$ be an $\I$-admissible filtration of ideals in $R.$ Then the following statements are equivalent.
\item[(1)] For all $\n\in\NN^s,$ $H_{\R_{++}}^1(\R(\fil))_{\n}=0.$
\item[(2)] For all $\ii,$ grade$(\gi)\geq 1.$ 
\bpf
Fix $i.$ Denote $\displaystyle\frac{\mathcal{R}'(\fil)}{\mathcal{R}'(\fil)(e_i)}$ by $G_i'(\fil).$ Consider the short exact sequence of $\R(\I)$-modules,
\beqn
0\longrightarrow \mathcal{R}'(\fil)(e_i)\longrightarrow \mathcal{R}'(\fil)\longrightarrow G_i'(\fil)\longrightarrow 0.
\eeqn This induces a long exact sequence of local cohomology modules, \beqnn\label{number}
 0\longrightarrow H_{{\R}_{++}}^0( G_i'(\fil))_{\n} \longrightarrow H_{{\R}_{++}}^1(\mathcal{R}'(\fil))_{\n+e_i} \longrightarrow\cdots.\eeqnn Since $H_{\R_{++}}^1(\R(\fil))_{\n}=0$ for all $\n\in\NN^s,$ by \cite[Proposition 4.2]{msv}, we have $H_{\R_{++}}^1(\R'(\fil))_{\n}=0$ for all $\n\in\NN^s.$ Hence from the exact sequence (\ref{number}) and \cite[Proposition 4.2]{msv}, for all $\n\in\NN^s,$ we get $$H_{{G_i(\fil)}_{++}}^0( G_i(\fil))_{\n}\cong H_{G_i(\fil)_{++}}^0( G_i'(\fil))_{\n}\cong H_{{\R}_{++}}^0( G_i'(\fil))_{\n}=0.$$
 \\Conversely, suppose grade$(\gi)\geq 1$ for all $\ii.$ By \cite[Theorem 3.3]{msv}, for all $\n\in\NN^s,$ $${\breve\fil(\n+e_i)\cap\fil(\n)}={\fil(\n+e_i)}.$$ We show that if $\breve\fil(\n)=\fil(\n)$ for some $\n\geq \underline 0$ then $\breve\fil(\m)=\fil(\m)$ for all $\m\geq\n.$ Let $t_i=m_i-n_i$ for all $\ii.$ For each $i,$ $$\fil(\n+e_i)=\breve\fil(\n+e_i)\cap\fil(\n)=\breve\fil(\n+e_i)\cap\breve\fil(\n)=\breve\fil(\n+e_i).$$ Continuing this process $t_i$ times for each $i$ we get $\breve\fil(\m)=\fil(\m).$ Since $\breve\fil(\underline 0)=\fil(\underline 0),$ by Proposition \ref{h1}, we get the required result. 
\eepf
\ep
\begin{remark}
Note that if $(R,\mm)$ is an analytically unramified local ring of dimension $d\geq 1,$ $\idl$ are $\mm$-primary ideals of $R$ then by \cite[Corollary 3.4]{msv}, for the filtration $\fil=\lbrace\overline{\I^{\n}}\rbrace_{\n\in\ZZ^s},$ $H_{\R_{++}}^1(\R(\fil))_{\n}=0$ for all $\n\in\NN^s.$
\end{remark}
\bp\label{bangalore}
Let $(R,\mm)$ be a \CM local ring of dimension $d\geq 2$ with infinite residue field and $\idl$ be $\mm$-primary ideals of $R.$ Let $\fil=\lbrace\fil(\n)\rbrace_{\n\in \ZZ^s}$ be an $\I$-admissible filtration of ideals in $R$ and $\mathcal {A}={\lbrace}{x_{ij}\in{I_i}: j=1,\ldots,d;\ii}\rbrace$ be a {{\gcrr}}of $\fil.$ Put $y_1=x_{11}\ldots x_{s1}.$ Let $H_{\R_{++}}^1(\R(\fil))_{\n}=0$ for all $\n\in\NN^s.$ Then \ben
{
\item[(1)] $y_{i1}=y_1+{\I}^{e+e_i}\in {G_{i}(\I)_{e}}$ is $\ggi$-regular for all $\ii.$  
\item[(2)] $(y_1)\cap\fil(\n)=y_1\fil(\n-e)$ for all $\n\geq e.$ 
}
\een
\bpf
$(1)$ Fix $i.$ Let $\m\geq e$ such that $(y_1)\cap\fil(\n)=y_1\fil(\n-e)$ for all $\n\geq \m.$  We show that $(\ggi)_{\n}\overset{.{y_{i1}}}\longrightarrow (\ggi)_{\n+e}$ is injective for all $\n\geq{\m}.$
Let $(z+{\fil}(\n+e_i)){y_{i1}}={\fil}(\n+e+e_i).$ Then $y_1z\in{{\fil}(\n+e+e_i)}.$ Since $\n\geq{\m},$  $z\in{{\fil}(\n+e_i)}.$ Hence by Propositions \ref{r} and \ref{grade}, $y_{i1}$ is a nonzerodivisor of $\ggi.$
\\$(2)$ For all $\ii,$ consider the Koszul complex $K_i.=K_i.({y_{i1}},\ggi):$ $$0\longrightarrow {\ggi}\longrightarrow {\ggi}(e)\longrightarrow 0.$$ The $\n$th component of this complex is $K_i.({y_{i1}},\ggi,\n):$ $$0\longrightarrow {\ggi}_{(\n)}\longrightarrow {\ggi}_{\n+e}\longrightarrow 0.$$ Hence for all $\ii$ and $\n\geq \underline 0,$ we have exact sequence of complexes $$0\longrightarrow {K_i.({y_{i1}},\ggi,\n)}\longrightarrow  C.(y_1,\fil(\n+e_i))\longrightarrow C.(y_1,\fil(\n))\longrightarrow 0,$$ which gives a long exact sequence of homology modules
\beqnn\label{sk}
\cdots\longrightarrow H_j(K_i.({y_{i1}},\ggi,\n))\longrightarrow  H_j(C.(y_1,\fil(\n+e_i)))\longrightarrow  H_j(C.(y_1,\fil(\n)))\longrightarrow\cdots.\eeqnn Since $y_{i1}$ is a nonzerodivisor of $\ggi$ for all $\ii,$ we have $H_1(K_i.({y_{i1}},\ggi,\n))=0$ for all $\n$ and $\ii.$ Since $H_1(C.(y_1,\fil(\underline 0)))=0,$ applying (\ref{sk}) several times for all $\ii,$ we get $H_1(C.(y_1,\n))=0$ for all $\n\geq \underline 0.$ Hence ${{{\fil}}(\n)}\cap{(y_1)}=y_{1}{{{\fil}}(\n-{e})}$ for all $\n\geq e.$
\eepf
\ep
\bp\label{m}
Let $(R,\mm)$ be a \CM local ring of dimension two and $\idl$ be $\mm$-primary ideals of $R.$ Let $\fil=\lbrace\fil(\n)\rbrace_{\n\in \ZZ^s}$ be an $\I$-admissible filtration of ideals in $R$ and $\mathcal {A}={\lbrace}{x_{ij}\in{I_i}: j=1,2;\ii}\rbrace$ be a \gcrr of $\fil.$ Let $y_j=x_{1j}\cdots x_{sj}$ for all $j=1,2$ and  $H_{\R_{++}}^1(\R(\fil))_{\n}=0$ for all $\n\in\NN^s.$ Then $H_2(C.(y_1,y_{2},\fil(\n)))=0$ for all $\n\geq \underline 0.$
\bpf
Fix $i.$ Let $y_{ij}=y_j+{\I}^{e+e_i}\in {G_{i}(\I)_{e}}$ for all $j=1,2$ and $\ii.$ Consider the Koszul complex ${K.(y_{i1},y_{i2},\ggi)}$:
$$0\longrightarrow \ggi\longrightarrow(\ggi(e))^2\longrightarrow\ggi(2e)\longrightarrow 0$$ whose $\n$th component is $$0\longrightarrow \ggi_{\n}\longrightarrow(\ggi_{\n+e})^2\longrightarrow\ggi_{\n+2e}\longrightarrow 0.$$
For all $\n\geq\underline{0},$ we have the exact sequence $$0\longrightarrow {K.(y_{i1},y_{i2},\ggi,\n)}\longrightarrow  C.(y_1,y_2,\fil(\n+e_i))\longrightarrow C.(y_1,y_2,\fil(\n))\longrightarrow 0.$$ This gives a long exact sequence of homology modules $$\cdots\longrightarrow {H_j(K.(y_{i1},y_{i2},\ggi,\n))}\longrightarrow {H_j(C.(y_1,y_2,\fil(\n+e_i)))}\longrightarrow {H_j(C.(y_1,y_2,\fil(\n)))}\longrightarrow\cdots$$ for all $\n\geq \underline 0.$ Since by Proposition \ref{bangalore}, $y_{i1}$ is a regular element in $\ggi$ for all $\ii,$ we have $H_2(K.({y_{i1}},{y_{i2}},\ggi))=0$ for all $\n$ and $\ii.$ Since $H_2(C.(y_1,y_2,\fil(\underline{0})))=0,$ using the above exact sequence several times for all $\ii,$ we get $H_2(C.(y_1,y_2,\fil(\n)))=0$ for all $\n\geq \underline 0.$
\eepf
\ep
\bt\label{vanishing}
Let $(R,\mm)$ be a \CM local ring of dimension $1\leq d\leq 2$ with infinite residue field, $\idl$ be $\mm$-primary ideals of $R$ and $\fil=\lbrace\fil(\n)\rbrace_{\n\in \ZZ^s}$ be an $\I$-admissible filtration of ideals in $R.$ Let $e_{(d-1)e_i}(\fil)=0$ for $\ii.$ Then 
\ben
\item For $d=1,$ $\po(\n)=\ho(\n)$ for all $\n\in\NN^s.$
\item For $d=2,$ if $H_{\R_{++}}^1(\R(\fil))_{\n}=0$ for all $\n\in\NN^s$ then $\po(\n)=\ho(\n)$ for all $\n\in\NN^s$ and $e_{\underline{0}}(\fil)=0.$\een
\bpf
$(1)$ For $d=1,$ since $e_{\underline{0}}(\fil)=0,$ we get $\po(\underline 0)=\ho(\underline 0).$ Hence by the difference formula \cite[Theorem 4.3]{msv}, we get $\lm_R\lf H_{\R_{++}}^1(\R(\fil))_{\underline{0}}\rg=P_{\fil}(\underline{0})-H_{\fil}(\underline{0})=0.$ Since $\dim R=1,$ by \cite[Lemma 2.11]{msv}, for all $\n\in\NN^s,$ $H_{\R_{++}}^1(\R(\fil))_{\n}=0.$ Therefore using the difference formula \cite[Theorem 4.3]{msv} again, we get $P_{\fil}(\n)-H_{\fil}(\n)=0$ for all $\n\in\NN^s.$ 
\\$(2)$ Let $d=2,$ $\mathcal {A}={\lbrace}{x_{ij}\in{I_i}: j=1,2;\ii}\rbrace$ be any \gcrr of $\fil,$ $y_j=x_{1j}\cdots x_{sj}$ for $j=1,2$ and $J=(y_1,y_2).$ Fix $i.$ Let $R'=R/(x_{i1})$ and $'$ denote the image of an ideal in $R'.$ For all large $\n,$ consider the following exact sequence $$
0\longrightarrow\frac{(\fil(\n):(x_{i1}))}{\fil(\n-e_i)}\longrightarrow\frac{R}{\fil(\n-e_i)}\overset{.x_{i1}}\longrightarrow\frac{R}{\fil(\n)}\longrightarrow\frac{R}{(x_{i1},\fil(\n))}\longrightarrow 0.$$ Since $\mathcal A$ is a good complete reduction, for all large $\n,$ $(\fil(\n):(x_{i1}))=\fil(\n-e_i)$ and
$$\lm\lf\frac{R}{(x_{i1},\fil(\n))}\rg=\lm\lf\frac{R}{\fil(\n)}\rg-\lm\lf\frac{R}{\fil(\n-e_i)}\rg.$$ Therefore for the filtration $\fil'=\{\fil(\n)R'\}_{\n\in\ZZ^s},$ we have $P_{\fil^\prime}(\n)=\po(\n)-\po(\n-e_i).$ This implies that the constant term of $P_{\fil^\prime}(\n)$ is $e_{e_i}(\fil)=0.$ 
\\Since ${\mathcal A}'=\lbrace{x_{i2}^\prime\in {I_i}: \ii }\rbrace$ is a complete reduction of $\fil',$ $J'=(y_2^\prime)$ and $\dim R'=1,$ by Theorem \ref{use}, Proposition \ref{a} and part $(1)$, we have $$\lm\lf\frac{\fil(\n+e)R'}{J'\fil(\n)R'}\rg=\Delta^1(P_{\fil^\prime}(\n)-H_{\fil^\prime}(\n))=0\mbox{ for all }\n\in\NN^s.$$ Thus we get $\fil(\n+e)=y_2\fil(\n)+((x_{i1})\cap\fil(\n+e))$ for all $\n\in\NN^s.$ We show that $(x_{i1})\cap\fil(\n+e)=x_{i1}\fil(\n+e-e_i)$ for all $\n\in\NN^s.$ It is clear that $x_{i1}\fil(\n+e-e_i)\subseteq (x_{i1})\cap\fil(\n+e).$ Let $ax_{i1}\in\fil(\n+e).$ Then $ay_1\in\fil(\n+2e-e_i).$ Since $H_{\R_{++}}^1(\R(\fil))_{\n}=0$ for all $\n\in\NN^s$ and $\mathcal A$ is a good complete reduction, by Proposition \ref{bangalore}, $a\in\fil(\n+e-e_i).$ Hence we get $$\fil(\n+e)=y_2\fil(\n)+x_{i1}\fil(\n+e-e_i)\mbox{ for all }\n\in\NN^s\mbox{ and }\ii.$$ We show that $\fil(\n+2e)=J\fil(\n+e)$ for all $\n\in\NN^s.$ Let $\n\in\NN^s.$ Then 
\beqn
\fil(\n+2e)&=& y_2\fil(\n+e)+x_{11}\fil(\n+2e-e_1)\\&=& y_2\fil(\n+e)+x_{11}(y_2\fil(\n+e-e_1)+x_{21}\fil(\n+2e-e_1-e_2))\\&\subseteq & y_2\fil(\n+e)+x_{11}x_{21}\fil(\n+2e-e_1-e_2)\\&\vdots&\\&\subseteq & y_2\fil(\n+e)+x_{11}\cdots x_{s1}\fil(\n+e)= (y_1,y_2)\fil(\n+e)= J\fil(\n+e).
\eeqn Since $H_{\R_{++}}^1(\R(\fil))_{\n}=0$ for all $\n\in\NN^s,$ by Theorem \ref{use} and Proposition \ref{m}, for all $\n\in\NN^s,$ we get \beqn
{\Delta}^2(\po(\n)-\ho(\n))=\lm\lf\frac{{\fil}(\n+2e)}{J{\fil}(\n+e)}\rg\eeqn and hence by Proposition \ref{a}, we have $\po(\n)=\ho(\n)$ for all $\n\in\NN^s.$ Putting $\n=\underline 0$ in the above equality we get $e_{\underline{0}}(\fil)=0.$
\eepf
\et
\bt
Let $(R,\mm)$ be a \CM local ring of dimension $1\leq d\leq 2$ with infinite residue field and $I,J$ be $\mm$-primary ideals of $R.$ Let $\fil=\lbrace\fil(r,s)\rbrace_{r,s\in \ZZ}$ be a $\ZZ^2$-graded $(I,J)$-admissible filtration of ideals in $R.$ Then the following statements are equivalent.\ben
\item $e_{(d-1)e_i}(\fil)=0$ for $i=1,2.$
\item $I$ and $J$ are generated by system of parameters, $\po(r,s)=\ho(r,s)$ for all $r,s\in\NN$ and $\fil(r,s)=I^rJ^s$ for all $r,s\in\ZZ.$ 
\item $e_{\alpha}(\fil)=0$ for $|\alpha|\leq d-1.$
\een
\bpf
$(1)\Rightarrow (2):$ Let $\mathcal F^{(1)}=\lbrace\fil(r,0)\rbrace_{r\in \ZZ}$ and $\mathcal F^{(2)}=\lbrace\fil(0,s)\rbrace_{s\in \ZZ}.$ Since $\fil$ is an $(I,J)$-admissible filtration, $\mathcal F^{(1)}$ and $\mathcal F^{(2)}$ are $I$-admissible and $J$-admissible filtrations respectively. By \cite{northcott}, \cite[Lemma 3.19]{marleythesis} and \cite[Theorem 5.5]{msv}, we have
$$0\leq e_1(I)\leq e_1(\mathcal F^{(1)})\leq e_{(d-1)e_1}(\fil)=0\mbox{ and }0\leq e_1(J)\leq e_1(\mathcal F^{(2)})\leq e_{(d-1)e_2}(\fil)=0.$$ Then by \cite[Theorem 3.21]{marleythesis}, we get $I$ and $J$ are generated by system of parameters, $\fil(r,0)=I^r$ and $\fil(0,s)=J^s$ for all $r,s\in\ZZ.$
\\Let $d=1.$ Then by Theorem \ref{vanishing}, $\po(r,s)=\ho(r,s)$ for all $r,s\in\NN.$ It is enough to prove that $\fil(r,s)=I^rJ^s$ for all $r,s\geq 1.$ Since $I,J$ are generated by system of parameters, for $r,s\geq 1,$ we have $$\lm\lf\frac{R}{\fil(r,s)}\rg=\po(r,s)=re(I)+se(J)=\lm\lf\frac{R}{I^r}\rg+\lm\lf\frac{R}{J^s}\rg=\lm\lf\frac{R}{I^r}\rg+\lm\lf\frac{I^r}{I^rJ^s}\rg=\lm\lf\frac{R}{I^rJ^s}\rg.$$ This implies $\fil(r,s)=I^rJ^s$ for all $r,s\geq 1.$
\\ Let $d=2.$ Since $I,J$ are parameter ideals, we have $$e(I)-e_{e_1}(\fil)=e(I)=\lm\lf\frac{R}{I}\rg=\lm\lf\frac{R}{\fil(1,0)}\rg\mbox{ and }e(J)-e_{e_2}(\fil)=e(J)=\lm\lf\frac{R}{J}\rg=\lm\lf\frac{R}{\fil(0,1)}\rg.$$ Therefore by \cite[Theorem 7.3]{msv}, we get $\po(r,s)=\ho(r,s)$ for all $r,s\in\NN.$ It is enough to prove that $\fil(r,s)=I^rJ^s$ for all $r,s\geq 1.$ By \cite[Theorem 7.3]{msv}, the joint reduction number of $\fil$ of type $e$ is zero. Let $(a,b)$ be a joint reduction of $\fil$ of type $e.$ Then $$\fil(r,s)=a\fil(r-1,s)+b\fil(r,s-1)\mbox{ for all }r,s\geq 1.$$  We use induction on $r+s.$ Let $r,s\geq 1.$ If $r+s=2$ then $r=s=1$ and $$\fil(1,1)=a\fil(0,1)+b\fil(1,0)=aJ+bI\subseteq IJ\subseteq \fil(1,1).$$ Let $r+s>2.$ Then $r\geq 2$ or $s\geq 2.$ Without loss of generality assume $r\geq 2.$ If $s=1$ then using induction we get $$\fil(r,1)=a\fil(r-1,1)+b\fil(r,0)=aI^{r-1}J+bI^r\subseteq I^rJ\subseteq \fil(r,1).$$ Hence we may assume $s\geq 2.$ Therefore $$\fil(r,s)=a\fil(r-1,s)+b\fil(r,s-1)=aI^{r-1}J^s+bI^rJ^{s-1}\subseteq I^rJ^s\subseteq \fil(r,s).$$
\\$(2)\Rightarrow (3):$ For $d=1,$ putting $r=s=0$ in the equation $\po(r,s)=\ho(r,s),$ we get the required result. Let $d=2.$ Since $\po(r,s)=\ho(r,s)$ for all $r,s\in\NN,$ we have $$e_{\underline 0}(\fil)=0,\mbox{ } e(I)-e_{e_1}(\fil)=\lm\lf\frac{R}{\fil(1,0)}\rg\mbox{ and }e(J)-e_{e_2}(\fil)=\lm\lf\frac{R}{\fil(0,1)}\rg.$$ Since $I,J$ are parameter ideals and $\fil(r,s)=I^rJ^s$ for all $r,s\in\ZZ,$ by \cite[Theorem 5.5]{msv}, we have $$e_{e_1}(\fil)=e_1(\mathcal F^{(1)})=e_1(I)=0\mbox{ and }e_{e_2}(\fil)=e_1(\mathcal F^{(2)})=e_1(J)=0.$$
\\$(3)\Rightarrow (1):$ It follows directly.

\eepf
\et
\bex
\rm{Let $R=k[|X,Y|].$ Then $R$ is a regular local ring of dimension two. Let $I=(X,Y^2)$ and $J=(X^2,Y).$ Then $I,J$ are complete parameter ideals in $R.$ Consider the filtration $\fil=\{\overline{I^rJ^s}\}_{r,s\in\ZZ}.$ Since $I,J$ are complete ideals, by \cite[Theorem 2$'$, Appendix 5]{zariski}, $I^r,J^s,I^rJ^s$ are complete ideals. By \cite[Theorem 1.2]{Rees}, $${e}_{e_1}(\fil)={\overline e}_1(I)=e_1(I)=0\mbox{ and }{e}_{e_2}(\fil)={\overline e}_1(J)=e_1(J)=0.$$ Since $${e}(I)-e_{e_1}(\fil)=e(I)=\lm\lf\frac{R}{I}\rg=\lm\lf\frac{R}{\overline I}\rg\mbox{ and }{e}(J)-e_{e_2}(\fil)=e(J)=\lm\lf\frac{R}{J}\rg=\lm\lf\frac{R}{\overline J}\rg,$$ by \cite[Theorem 7.3]{msv}, we get $\po(r,s)=\ho(r,s)$ for all $r,s\in\NN.$
}
\eex
\section{Postulation and reduction vectors in dimension one}
Let $(R,\mm)$ be a \CM local ring of dimension one with infinite residue field and $\idl$ be $\mm$-primary ideals of $R.$ Let $\fil=\lbrace\fil(\n)\rbrace_{\n\in \ZZ^s}$ be an $\I$-admissible filtration of ideals in $R.$ In this section we prove that the set of reduction vectors of $\fil$ with respect to any complete reduction is the same as the set of postulation vectors of $\fil.$ Thus the set of reduction vectors of $\fil$ with respect to any \crr is independent of the choice of complete reduction. Then we show that the complete reduction number of $\fil$ \wrt any complete reduction is independent of choice of complete reduction. 
\bt\label{dimension one}
Let $(R,\mm)$ be a \CM local ring of dimension one with infinite residue field and $\idl$ be $\mm$-primary ideals of $R.$ Let $\fil=\lbrace\fil(\n)\rbrace_{\n\in \ZZ^s}$ be an $\I$-admissible filtration of ideals in $R$ and $\mathcal {A}=\lbrace{a_i\in{I_i}:\ii}\rbrace$ be a complete reduction of $\fil.$ Then $$\PP(\fil)\subseteq\NN^s\mbox{  and  }\PP(\fil)={\mathcal{R}}_{\mathcal A}(\fil).$$ Moreover, the set $\R_{\mathcal A}(\fil)$ is independent of any complete reduction $\mathcal A$ of $\fil.$
\bpf 
First we prove that $\PP(\fil)\subseteq\NN^s.$ Suppose there exists $\n\in\ZZ^s\setminus\NN^s$ such that $\n\in\PP(\fil).$ Then there exists at least one $i\in\{1,\ldots,s\}$ such that $n_i<0.$ Therefore $$\po(\n+e_i)=\lm\lf\frac{R}{\fil(\n+e_i)}\rg=\lm\lf\frac{R}{\fil(\n)}\rg=\po(\n)$$ implies $e_0(I_i)=\po(\n+e_i)-\po(\n)=0.$ This contradicts to the fact that $e_0(I_i)>0.$ Thus $\PP(\fil)\subseteq\NN^s.$ 
\\Let $J=(a_{1}\cdots a_{s}).$ By Theorem \ref{use}, for all $\n\geq \underline{0},$ \beqn\label{p1}
{\Delta}^1(\po(\n)-\ho(\n))=\lm\lf\frac{{\fil}(\n+e)}{J{\fil}(\n)}\rg. \eeqn Hence by Proposition \ref{a}, we get the required result.
\eepf
\et
\bex\rm{
Let $R=k[|t^3,t^4,t^5|].$ Then $R$ is a one-dimensional \CM local ring with unique maximal ideal $\mm=(t^3,t^4,t^5)$. Consider $I=(t^3,t^4)$ and $J=(t^3).$ Then $JI^2=I^3.$ Since $(t^6)(IJ)^2=(IJ)^3,$ $\mathcal A=\begin{pmatrix}
  t^3\\
  t^3
 \end{pmatrix}$ is a \crr for the filtration $\I=\lbrace I^rJ^s\rbrace_{r,s\in\ZZ}.$ We have $\lm\lf{\frac{R}{J^n}}\rg=3n$ for all $n\in\NN.$ Now $\lm\lf{\frac{R}{I}}\rg=2,$ $\lm\lf{\frac{R}{I^2}}\rg=4$ and for $n\geq 3,$ 
$$
I^n =J^{n-2}I^2=(t^{3n-6})(t^6,t^7,t^8)=(t^{3n},t^{3n+1},t^{3n+2}).$$
 Hence for all $n\geq 2,$ $\lm\lf{\frac{R}{I^n}}\rg=3n-2.$ Let $P_I(n)=ne_0(I)-e_1(I),$ $P_J(n)=ne_0(J)-e_1(J),$ $P_{IJ}(n)=ne_0(IJ)-e_1(IJ)$ and $\pf(\n)=n_1e_0(I)+n_2e_0(J)-{e_{\underline 0}(\I)}$ denote the Hilbert polynomials of $I,$ $J,$ $IJ$ and $\I$ respectively where $n\in\ZZ$ and $\n=(n_1,n_2)\in{\ZZ}^2.$ Then $e_0(I)=e_0(J)=3,$ $e_1(I)=2$ and $e_1(J)=0.$ Now by Lemma \ref{one}, ${e_{\underline 0}(\I)}=e_1(IJ).$ For large $n,$
$$
P_{IJ}(n)= ne_0(IJ)-e_1(IJ)=\lm\lf{\frac{R}{(IJ)^n}}\rg = \lm\lf{\frac{R}{I^{2n}}}\rg= 6n-2.$$
 Hence ${e_{\underline 0}(\I)}=e_1(IJ)=2.$ This implies $\pf(\n)=3n_1+3n_2-2.$
\\ Since $(t^6)(IJ)=IJ^3\neq I^2J^2$ and $(t^6)(IJ)^2=I^2J^4=I^3J^3,$ we have ${r}_{\mathcal A}(\I)=2.$
\\ Note that $(t^6)I^2=J^2I^2=JI^3$ and $(t^6)I=J^2I\neq I^2J.$ Let $(1,n)\in{\ZZ}^2$ such that $n\geq 1.$ Then $(t^6)IJ^n=IJ^{n+2}=(t^3,t^4)(t^{3n+6})=(t^{3n+9},t^{3n+10})\neq (t^{3n+9},t^{3n+10},t^{3n+11})=I^2J^{n+1}.$ Hence ${\R}_{\mathcal A}(\I)=\lbrace{\m\in\NN^2\mid\m\geq (2,0)}\rbrace.$
\\ For all $n_1\geq 2$ and $n_2\geq 0,$ $$\pf(n_1,n_2)=3n_1+3n_2-2=\lm\lf{\frac{R}{I^{n_1+n_2}}}\rg=\lm\lf{\frac{R}{I^{n_1}J^{n_2}}}\rg=\hf(n_1,n_2)$$ and $\pf(1,0)=1\neq 2=\lm\lf{\frac{R}{I}}\rg.$ Let $(1,n)\in{\ZZ}^2$ such that $n\geq 1.$ Then $\pf(1,n)=3n+1$ and $\hf(1,n)=\lm\lf{\frac{R}{IJ^n}}\rg=\lm\lf{\frac{R}{(t^{3n+3},t^{3n+4})}}\rg=3n+2.$ Hence for all $\n=(1,n)\mbox{ where }n\geq0,$ $\pf(1,n)\neq \hf(1,n).$ Thus ${\R}_{\mathcal A}(\I)=\PP(\I).$ }
\eex
In the following example we show that we cannot drop the condition of Cohen-Macaulayness in Theorem \ref{dimension one}.
\bex\rm{
Let $R=\displaystyle\frac{k[|X,Y|]}{(X^2,XY)}.$ Then $R$ is a one-dimensional local ring which is not Cohen-Macaulay. Consider the ideals $I=(x,y)$ and $J=(y)$ of $R.$ Then $JI=I^2.$ Since $(y^2)(IJ)=I^2J^2,$ $\mathcal A=\begin{pmatrix}
  y\\
  y
 \end{pmatrix}$ is a \crr for the filtration $\I=\lbrace I^rJ^s\rbrace_{r,s\in\ZZ}.$ We have $\lm\lf{\frac{R}{J^n}}\rg=n+1$ for all $n\geq 1.$
  Now $\lm\lf{\frac{R}{I}}\rg=1$ and for $n\geq 2,$ $I^n = J^{n-1}I=(y^{n-1})(x,y)=(y^{n}).$ Hence for all $n\geq 2,$ $\lm\lf{\frac{R}{I^n}}\rg=n+1.$ Let $P_I(n)=ne_0(I)-e_1(I),$ $P_J(n)=ne_0(J)-e_1(J),$ $P_{IJ}(n)=ne_0(IJ)-e_1(IJ)$ and $\pf(\n)=n_1e_0(I)+n_2e_0(J)-{e_{\underline 0}(\I)}$ denote the Hilbert polynomials \wrt $I,$ $J,$ $IJ$ and $\I$ respectively where $n\in\ZZ$ and $\n=(n_1,n_2)\in{\ZZ}^2.$ Then $e_0(I)=e_0(J)=1.$ Now by Lemma \ref{one}, ${e_{\underline 0}(\I)}=e_1(IJ).$ For large $n,$
$P_{IJ}(n)= ne_0(IJ)-e_1(IJ)=\lm\lf{\frac{R}{(IJ)^n}}\rg= \lm\lf{\frac{R}{I^{2n}}}\rg= 2n+1.$ Hence ${e_{\underline 0}(\I)}=e_1(IJ)=-1.$ This implies $\pf(\n)=n_1+n_2+1.$ Now $IJ=(x,y)(y)=(y^2).$ Hence ${r}_{\mathcal A}(\I)=0.$ This implies ${\R}_{\mathcal A}(\I)=\NN^2.$ But $\pf(0,0)=1\neq 0=\hf(0,0).$ Hence ${\R}_{\mathcal A}(\I)\neq\PP(\I).$}
\eex
\bt\label{eye}
Let $(R,\mm)$ be a one-dimensional \CM local ring and $\idl$ be $\mm$-primary ideals of $R.$ Let $\fil=\lbrace\fil(\n)\rbrace_{\n\in \ZZ^s}$ be an $\I$-admissible filtration of ideals in $R.$ Then the complete reduction number of $\fil$ \wrt any \crr is independent of choice of complete reduction of $\fil$.
\bpf
Let $\mathcal {A}=\lbrace{a_i\in{I_i}:\ii}\rbrace$ be a \crr of $\fil,$ $J=(a_{1}\cdots a_{s})$ and ${r}_{\mathcal A}(\fil)=k.$ First we show that $$k=\min{\lbrace}{\max{\lbrace}{t_1},\ldots,{t_s}:{\underline{t}}=({t_1},\ldots,{t_s})\in{{\R}_{\mathcal A}(\fil)}{\rbrace}}{\rbrace}.$$
If $k=0$ then it is true. Suppose $k\geq 1.$ Let $\n\in{\NN}^s$ such that $n_i< k$ for all $\ii$ and $\n\in {{\R}_{\mathcal A}(\fil)}.$ Let $u=\max{\lbrace}{n_1},\ldots,{n_s}{\rbrace}.$ Then $u<k$ and ${\n}\leq ue\leq (k-1)e.$ Hence $J\fil(\m)=\fil(\m+e)$ for all $\m\geq (k-1)e.$ This contradicts to the fact that $k$ is the complete reduction number of $\fil$ \wrt $\mathcal A.$ Thus ${\underline{t}}\in{{\R}_{\mathcal A}(\fil)}$ implies $t_i\geq k$ for at least one $i.$ Since $J{\fil}(\n)={\fil}(\n+e)$ for all $\n\geq ke,$ there exists ${\rrr}\in{{\R}_{\mathcal A}(\fil)}$ such that $\max{\lbrace}{r_1},\ldots,{r_s}{\rbrace}=k.$ Thus $k=\min{\lbrace}{\max{\lbrace}{t_1},\ldots,{t_s}:{\underline{t}}=({t_1},\ldots,{t_s})\in{{\R}_{\mathcal A}(\fil)}{\rbrace}}{\rbrace}.$
 Hence by Theorem \ref{dimension one}, we get the required result. 
\eepf
\et
\section{Postulation and reduction vectors in dimension two}
Let $(R,\mm)$ be a \CM local ring of dimension two with infinite residue field and $\idl$ be $\mm$-primary ideals of $R.$ Let $\fil=\lbrace\fil(\n)\rbrace_{\n\in \ZZ^s}$ be an $\I$-admissible filtration of ideals in $R.$ In this section we provide a relation between the reduction vectors of $\fil$ with respect to any good complete reduction and the postulation vectors of $\fil.$ For a bigraded filtration $\fil,$ we prove a result which relates the Cohen-Macaulayness of the bigraded Rees algebra, the complete reduction number, reduction numbers and the joint reduction number.
\bt\label{j}
Let $(R,\mm)$ be a \CM local ring of dimension two with infinite residue field and $\idl$ be $\mm$-primary ideals of $R$ and $s\geq 2.$ Let $\fil=\lbrace\fil(\n)\rbrace_{\n\in \ZZ^s}$ be an $\I$-admissible filtration of ideals in $R$ and $\mathcal {A}={\lbrace}{x_{ij}\in{I_i}: j=1,2;\ii}\rbrace$ be a \gcrr of $\fil.$ Let $H_{\R_{++}}^1(\R(\fil))_{\n}=0$ for all $\n\geq\underline 0.$ Then $\PP(\fil)\subseteq\NN^s$ and there exists a one-to-one correspondence 
$$f:{\displaystyle{\Large{\Large{\PP(\fil)\longleftrightarrow\lbrace{\rrr\in{\R_{\mathcal A}(\fil)}\mid\rrr\geq e}\rbrace}}}}$$ defined by $f(\n)=\n+e$ where $f^{-1}(\rrr)=\rrr-e.$ 
\bpf 
First we prove that $\PP(\fil)\subseteq\NN^s.$ Suppose there exists $\n\in\ZZ^s\setminus\NN^s$ such that $\n\in\PP(\fil).$ Then there exists at least one $i\in\{1,\ldots,s\}$ such that $n_i<0.$ Therefore for any $j\in\{1,\ldots,s\}$ with $j\neq i$ and $l\geq 0,$  $$\po(\n+le_j+e_i)=\lm\lf\frac{R}{\fil(\n+le_j+e_i)}\rg=\lm\lf\frac{R}{\fil(\n+le_j)}\rg=\po(\n+le_j).$$ Thus for $l=1,$ we get \beqnn\label{ramanujan} 0&=&\po(\n+e_j+e_i)-\po(\n+e_j)\nonumber\\&=&(n_i+1)e_0(I_i)+\sum\limits_{k\neq i,j}n_ke_{e_k+e_i}(\fil)+(n_j+1)e_{e_j+e_i}(\fil)-e_{e_i}(\fil).\eeqnn Then for $l=0,$ by equation (\ref{ramanujan}), we get \beqn 0 &=&\po(\n+e_i)-\po(\n)\\&=&(n_i+1)e_0(I_i)+\sum\limits_{k\neq i}n_ke_{e_k+e_i}(\fil)-e_{e_i}(\fil)\\&=&-e_{e_j+e_i}(\fil).\eeqn This contradicts to the fact that $e_{e_j+e_i}(\fil)>0$ \cite[Theorem 2.4]{rees3}. Hence $\PP(\fil)\subseteq\NN^s.$
\\Let $y_j=x_{1j}\cdots x_{sj}$ for $j=1,2$ and $J=(y_1,y_2).$ Then by Theorem \ref{use} and Proposition \ref{m}, for all $\n\geq \underline{0},$ \beqn
{\Delta}^2(\po(\n)-\ho(\n))=\lm\lf\frac{{\fil}(\n+2e)}{J{\fil}(\n+e)}\rg.\eeqn Hence by Proposition \ref{a}, we get the required result.
\eepf
\et
\bt\label{b}
Let $(R,\mm)$ be a \CM local ring of dimension two with infinite residue field and $\idl$ be $\mm$-primary ideals of $R$ and $s\geq 2.$ Let $\fil=\lbrace\fil(\n)\rbrace_{\n\in \ZZ^s}$ be an $\I$-admissible filtration of ideals in $R$ and $H_{\R_{++}}^1(\R(\fil))_{\n}=0$ for all $\n\geq\underline 0.$ Then the following statements are equivalent.
\ben
{
\item[(1)] $\PP(\fil)=\NN^s,$ i.e. $\po(\n)=\ho(\n)$ for all $\n\geq\underline 0.$
\item[(2)] $r_{\mathcal A}(\fil)\leq {1}$ for any \gcrr $\mathcal {A}$ of $\fil.$
\item[$(2')$] There exists a \gcrr $\mathcal {A}$ of $\fil$ such that $r_{\mathcal A}(\fil)\leq {1}.$
}
\een
\bpf 
$(1)\Rightarrow(2):$
Let $\po(\n)=\ho(\n)$ for all $\n\geq\underline{0}$ and $\mathcal {A}={\lbrace}{x_{ij}\in{I_i}: j=1,2;\ii}\rbrace$ be any \gcrr of $\fil.$ Let $y_j=x_{1j}\cdots x_{sj}$ for $j=1,2$ and $J=(y_1,y_2).$ Then by Theorem \ref{use}, Propositions \ref{m} and \ref{a}, for all $\n\geq\underline{0},$
$$\lm\lf\frac{{\fil}(\n+2e)}{J{\fil}(\n+e)}\rg=0\Longrightarrow J{\fil}(\n+e)={\fil}(\n+2e).$$ Hence $r_{\mathcal A}(\fil)\leq {1}.$\\ The implication $(2)\Rightarrow(2')$ is trivial. 
\\$(2')\Rightarrow(1):$  Suppose there exists a \gcrr $\mathcal {A}={\lbrace}{x_{ij}\in{I_i}: j=1,2;\ii}\rbrace$ of $\fil$ such that $r_{\mathcal A}(\fil)\leq {1}.$ Let $y_j=x_{1j}\cdots x_{sj}$ for $j=1,2$ and $J=(y_1,y_2).$ Then again by Theorem \ref{use} and Proposition \ref{m}, for all $\n\geq\underline{0},$  ${\Delta}^2(\po(\n)-\ho(\n))=0.$ Now using Proposition \ref{a}, we get $\po(\n)-\ho(\n)=0$ for all $\n\geq\underline{0}.$ Since $\PP(\fil)\subseteq\NN^s,$ we get the required result.
\eepf
\et
In the following example we show that we cannot drop the condition on $H_{\R_{++}}^1(\R(\fil))_{\n}$ in Theorem \ref{b}.
\bex\rm{
Let $R=k[|X,Y|].$ Then $R$ is a two-dimensional \CM local ring with unique maximal ideal $\mm=(X,Y).$ Let $I=\mm^2$ and $J=(X^2,Y^2).$ Since $(X^{4},Y^{4})IJ=(X^{4},Y^{4})\mm^{4}=\mm^{8}=I^2J^2,$ we have $\mathcal A=\begin{pmatrix}
  X^{2} & Y^{2}\\
  X^{2} & Y^{2}
 \end{pmatrix}$ is a \crr for the filtration $\I=\lbrace I^rJ^s\rbrace_{r,s\in\ZZ}.$ By \cite[Proposition 1.2.2]{HH}, for all large $\n=(n_1,n_2)\in\NN^2,$ we get $$(X^4)\cap I^{n_1}J^{n_2}=X^4(I^{n_1}J^{n_2}:(X^4))=X^4I^{n_1-1}J^{n_2-1}.$$
 Hence $\mathcal A$ is a \gcrr for the filtration $\I.$ Note that, since $\widebreve{{\I}^{e_2}}=\breve{J}=(I^k{J}^{1+k}:{I^kJ^k})$ for some large $k,$ $JI=I^2$ and $\mm$ is parameter ideal, we have $$\breve{J}=(I^{2k+1}:I^{2k})=(\mm^{4k+2}:\mm^{4k})=\mm^2\neq J\mbox{ and hence }H_{\R_{++}}^1(\R(\I))_{(0,1)}\neq 0.$$
 Since $(X^{4},Y^{4})IJ=I^2J^2,$ we get ${r}_{\mathcal A}(\I)\leq 1.$
 \\ For $n\geq 1,$ $$\lm\lf\frac{R}{I^n}\rg=\lm\lf\frac{R}{\mm^{2n}}\rg=\binom{2n+1}{2}=4\binom{n+1}{2}-n=P_I(n),$$ Since $J$ is parameter ideal and $JI=I^2,$ $$\lm\lf\frac{R}{J^n}\rg=4\binom{n+1}{2}=P_{J}(n),$$ $$\lm\lf\frac{R}{(IJ)^n}\rg=\lm\lf\frac{R}{\mm^{4n}}\rg=\binom{4n+1}{2}=16\binom{n+1}{2}-6n=P_{IJ}(n)$$ and
$$\lm\lf\frac{R}{I^{2n}J^{n}}\rg=\lm\lf\frac{R}{\mm^{6n}}\rg=\binom{6n+1}{2}=36\binom{n+1}{2}-15n=P_{I^2J}(n).$$ Hence $e_0(I)=e_0(J)=4.$ Now for large $n,$ $P_{IJ}(n)=\lm\lf\frac{R}{(IJ)^n}\rg=P_{\I}(ne)$ and $P_{I^2J}(n)=\lm\lf\frac{R}{I^{2n}J^{n}}\rg=P_{\I}(2n,n).$ Comparing the coefficients on both sides we get 
$$\pf(r,s)=4\binom{r+1}{2}+4\binom{s+1}{2}+4rs-r-s.$$ Then $\pf(0,1)=3\neq 4=\lm\lf\frac{R}{J}\rg.$ 
}\eex
\begin{theorem}\label{un}
Let $(R,\mm)$ be a \CM local ring of dimension two with infinite residue field and $I,J$ be $\mm$-primary ideals of $R.$ Let $\fil=\lbrace\fil(\n)\rbrace_{\n\in \ZZ^2}$ be a $\ZZ^2$-graded $(I,J)$-admissible filtration of ideals in $R.$ Then the following statements are equivalent.
\ben
{\item[(1)] The Rees algebra $\R(\fil)$ is Cohen-Macaulay.
\item[(2)] $\PP(\fil)=\NN^2,$ i.e. $\po(\n)=\ho(\n)$ for all $\n\geq\underline 0.$
\item[(3)] For any \gcrr $\mathcal {A}$ of $\fil,$ $r_{\mathcal A}(\fil)\leq {1}$ and $H_{\R_{++}}^1(\R(\fil))_{\n}=0$ for all $\n\geq\underline 0.$
\item[$(3')$] There exists a \gcrr $\mathcal {A}$ of $\fil$ such that $r_{\mathcal A}(\fil)\leq {1}$ and $H_{\R_{++}}^1(\R(\fil))_{\n}=0$ for all $\n\geq\underline 0.$
\item[(4)] For the filtrations $\mathcal F^{(i)}=\lbrace\fil(ne_i)\rbrace_{n\in \ZZ},$ $r(\mathcal F^{(i)}) \leq 1$ where $i=1,2$ and the joint reduction number of $\mathcal F$ of type $e$ is zero. 
}
\een
\bpf
The implications $(1)\Rightarrow(2)\Rightarrow(4)\Rightarrow (1)$ and $(3)\Rightarrow(3')\Rightarrow (2)$ follow from \cite[Theorem 7.3]{msv} and Theorem \ref{b} respectively. Enough to show $(1)\Rightarrow(3).$ Suppose $\mathcal R(\mathcal F)$ is Cohen-Macaulay. Then by \cite[Proposition 7.2,Theorem 2.15]{msv}, we get $H_{\R_{++}}^1(\R(\fil))_{\n}=0$ for all $\n\geq\underline 0.$ Since $(1)$ implies $(2),$ by Theorem \ref{b}, we get the required result. 
\eepf
\end{theorem}
The following example illustrates Theorems \ref{j}, \ref{b} and \ref{un}.
\bex\rm{
Let $R=K[|X,Y|].$ Then $R$ is a 2-dimensional \CM local ring with unique maximal ideal $\mm=(X,Y).$ Let $I=\mm^2$ and $J=\mm^3.$  Since $(X^{5},Y^{5})IJ=(X^{5},Y^{5})\mm^{5}=\mm^{10}=I^2J^2,$ we have $\mathcal A=\begin{pmatrix}
  X^{2} & Y^{2}\\
  X^{3} & Y^{3}
 \end{pmatrix}$ is a \crr for the filtration $\I=\lbrace I^rJ^s\rbrace_{r,s\in\ZZ}$ and ${r_{\mathcal A}(\I)}\leq 1.$ By \cite[Proposition 1.2.2]{HH}, for all large $\n=(n_1,n_2)\in\NN^2,$ we get $$(X^{5})\cap I^{n_1}J^{n_2}=X^{5}(I^{n_1}J^{n_2}:(X^{5}))=X^5I^{n_1-1}J^{n_2-1}.$$ Hence $\mathcal A$ is a \gcrr for the filtration $\I.$
 \\We show that $H_{\R_{++}}^1(\R(\I))_{\n}=0$ for all $\n\geq\underline 0.$  For $\n=(n_1,n_2)\geq\underline{0},$ we have $$\widebreve{{\I}^{\n}}=({\I}^{\n+ke}:{\I}^{ke})=(\mm^{2n_1+3n_2+5k}:\mm^{5k})$$ for some large $k.$ Since $\mm$ is parameter ideal, $\widebreve{{\I}^{\n}}=\mm^{2n_1+3n_2}={\I}^{\n}$ for all $\n\geq\underline 0.$ 
 \\ For $n\geq 1,$ $$\lm\lf\frac{R}{I^n}\rg=\lm\lf\frac{R}{\mm^{2n}}\rg=\binom{2n+1}{2}=4\binom{n+1}{2}-n=P_I(n),$$ $$\lm\lf\frac{R}{J^n}\rg=\lm\lf\frac{R}{\mm^{3n}}\rg=\binom{3n+1}{2}=9\binom{n+1}{2}-3n=P_{J}(n),$$ $$\lm\lf\frac{R}{(IJ)^n}\rg=\lm\lf\frac{R}{\mm^{5n}}\rg=\binom{5n+1}{2}=25\binom{n+1}{2}-10n=P_{IJ}(n)$$ and
$$\lm\lf\frac{R}{I^nJ^{2n}}\rg=\lm\lf\frac{R}{\mm^{8n}}\rg=\binom{8n+1}{2}=64\binom{n+1}{2}-28n=P_{IJ^2}(n).$$  
   Hence $e_0(I)=4$ and $e_0(J)=9.$ Now for large $n,$ $P_{IJ}(n)=\lm\lf\frac{R}{(IJ)^n}\rg=P_{\I}(ne)$ and $P_{IJ^2}(n)=\lm\lf\frac{R}{I^nJ^{2n}}\rg=P_{\I}(n,2n).$ Comparing the coefficients on both sides we get \beqnn\label{polynomial}\pf(n_1,n_2)=4\binom{n_1+1}{2}+9\binom{n_2+1}{2}+6n_1n_2-n_1-3n_2.\eeqnn Hence by \cite[Theorem 6.2]{msv}, the joint reduction number of $\I$ of type $e$ is zero. Since $(X^{2},Y^{2})I=I^2$ and $(X^{3},Y^{3})J=J^2,$ we have $r(\I^{(i)}) \leq 1$ for $i=1,2,$ where $\I^{(1)}=\{I^r\}_{r\in\ZZ}$ and $\I^{(2)}=\{J^s\}_{s\in\ZZ}.$  By \cite[Corollary 2.3, Corollary 2.4]{hhr}, $\R(\I^{(1)})$ and $\R(\I^{(2)})$ are Cohen-Macaulay. Hence by \cite[Theorem 7.1]{msv}, $\R(\I)$ is Cohen-Macaulay.
\\Using (\ref{polynomial}), we get $\pf(n_1,n_2)=\hf(n_1,n_2)$ for all $\n=(n_1,n_2)\in\NN^2.$ Thus $\mathcal P(\I)=\NN^2.$ Since $(X^{5},Y^{5})IJ=I^2J^2,$ $(X^{5},Y^{5})I\neq I^2J,$ $(X^{5},Y^{5})J\neq IJ^2,$ $(X^{5},Y^{5})I^2= I^3J$ and $(X^{5},Y^{5})J^2= IJ^3,$  we have $\R_{\mathcal A}(\I)=\{\m\in\NN^2\mid\m\geq e\}\cup\{\m\in\NN^2\mid\m\geq 2e_1\}\cup\{\m\in\NN^2\mid\m\geq 2e_2\}.$
 
}\eex

\end{document}